\documentclass[11pt]{article}
\usepackage{graphicx}
\usepackage{geometry}

\usepackage[T1]{fontenc}
\usepackage[latin1]{inputenc}

\usepackage{color}

\def \bop {\noindent\textbf{Proof. }}
\def \eop {\hbox{}\nobreak\hfill
\vrule width 2mm height 2mm depth 0mm
\par \goodbreak \smallskip}
\newcommand{\1}{1\!\!1}

\newcommand{\eps}  {\varepsilon}

\usepackage{amsmath}
\usepackage{amssymb}
\usepackage{amsfonts}
\usepackage{amsthm,amscd}

\newtheorem{theorem}{Theorem}[section]
\newtheorem{lemma}[theorem]{Lemma}
\newtheorem{corollary}[theorem]{Corollary}
\newtheorem{proposition}[theorem]{Proposition}
\theoremstyle{definition}

\newtheorem{remark}{Remark}[section]

\numberwithin{equation}{section}
\def\N{\mathbb{N}}
\def\R{\mathbb{R}}
\def\P{\mathbb{P}}
\def\E{\mathbb{E}}

\begin{document}
\title{ { Averaging for SDE-BSDE with null recurrent fast component Application
to homogenization in a non periodic media }}

\author {\quad K. Bahlali$^1$,\, A. Elouaflin$^2$,\,  E.
Pardoux$^3$}
 \maketitle


\begin{center}
{$^{1}$  Université de Toulon, B.P. 20132, 83957 La Garde Cedex,
France
 \\ \small{e-mail:
bahlali@univ-tln.fr}}
 \\
{$^2$ UFRMI, Universit\'e d'Abidjan,  22 BP 582 Abidjan, C\^{o}te
d'Ivoire \\
\small{e-mail:  elabouo@yahoo.fr}}
\\
{$^3$ Aix-Marseille Universit\'e, CNRS, Centrale Marseille,\\ I2M, UMR 7373 13453 Marseille, France \\
\small{e-mail: etienne.pardoux@univ-amu.fr}}
\end{center}



\begin{abstract}

We establish an averaging principle for a family of solutions
$(X^{\varepsilon}, Y^{\varepsilon})$ $ :=$ $(X^{1,\,\varepsilon},\,
X^{2,\,\varepsilon},\, Y^{\varepsilon})$ of a system of SDE-BSDE
with a null recurrent fast component $X^{1,\,\varepsilon}$. In
contrast to the classical periodic case, we can not rely on an
invariant probability and the slow forward component
$X^{2,\,\varepsilon}$ cannot be approximated by a diffusion process.
On the other hand, we assume that the coefficients admit a limit in a
\`{C}esaro sense. In such a case, the limit coefficients may have
discontinuity. We show that we can approximate the triplet
$(X^{1,\,\varepsilon},\,  X^{2,\,\varepsilon},\, Y^{\varepsilon})$ by
a system of SDE-BSDE $(X^1, X^2, Y)$ where $X := (X^1, X^2)$ is a
Markov diffusion which is the unique (in law) weak  solution of the
averaged forward  component and $Y$ is the unique solution to the
 averaged backward component. This is done with a backward component whose
generator depends on the variable $z$. As
application, we establish an homogenization result for semilinear
PDEs when the coefficients can be neither periodic nor ergodic.  We
show that the averaged BDSE is related to the averaged PDE via a
probabilistic representation of the (unique) Sobolev $
W_{d+1,\text{loc}}^{1,2}(\R_+\times\R^d)$--solution of the limit
PDEs. Our approach combines PDE methods and probabilistic arguments
which are based on stability property and weak convergence of BSDEs
in the S-topology.
\end{abstract}

\noindent {\bf Keys words}: \textit{SDE, BSDEs and PDES with
discontinuous coefficients, weak convergence of SDEs and BSDEs,
homogenization, S-topology, Averaging in \`{C}esaro sence, Sobolev
Spaces, Sobolev
solution to semilinear PDEs.}

\section{Introduction}

The averaging of stochastic differential equations (SDE)
as well as the homogenization of a partial differential equation
(PDE) is a process which consists in showing the convergence of the solution of an equation with rapidly varying
coefficients towards an equation with simpler (e.g. constant) coefficients.

 The two classical situations which were
mainly studied are the cases of deterministic periodic and random
stationary coefficients. These two situations are based on the existence of an
 invariant probability measure for some underlying process. The averaged coefficients are
then determined as a certain "means" with respect to this invariant
probability measure.

There is a vast literature on the homogenization of PDEs with
periodic coefficients, see for example the monographs \cite{BLP, FR,
PAN} and the references therein.
 There also exist numerous works on
averaging of stochastic differential equations with periodic
structures and its connection with homogenization of second order
partial differential equations (PDEs). Closer to our concern here,
 we can quote in particular
 \cite{BI, BHP, BH, D, EO, I,
L, EP2, P} and the references therein.

In contrast to these two classical situations (deterministic
periodic and random stationary coefficients) which were mainly
studied, we consider in this paper a different situation, building
upon earlier results of \cite{KK} and more recently those of
\cite{BEP, BDES}.  We extend the results of \cite{KK} to systems of
SDE-BSDEs and those of \cite{BEP, BDES} to the case where the
generator $f$ of the BSDE component depends upon the second unknown of the BSDE. As a consequence, we derive an homogenization result for
semilinear PDEs when the nonlinear part depends  on the solution as
well as on its gradient.

In \cite{KK}, Khasminskii \& Krylov consider the averaging of the
following family of diffusions process indexed by $\varepsilon$,

\begin{equation}\label{SDEH01}
 \left\{
\begin{aligned}
X^{1,x,\varepsilon}_s&=x_1+\int_0^s\varphi(\frac{X^{1,x,\varepsilon}_r}
{\varepsilon},\,X^{2,x,\varepsilon}_r)dW_r,\\
X^{2,x,\varepsilon}_s&=x_2+\int_0^s \tilde b
(\frac{X^{1,x,\varepsilon}_r}{\varepsilon},\,X^{2,x,\varepsilon}_r)dr
+\int_0^t\tilde\sigma (\frac{X^{1,x,\varepsilon}_r}
{\varepsilon},\,X^{2,x,\varepsilon}_r)d {W}_r,
\end{aligned}
\right.
\end{equation}
where $X^{1,x,\varepsilon}$ is a one-dimensional null-recurrent fast
component and $X^{2,x,\varepsilon}_t$  is a $d$--dimensional slow
component. The function $\varphi = (\varphi_1,..., \varphi_k)$ [resp. $\tilde\sigma = (\tilde\sigma_{ij})_{i,j}$, resp.
$\tilde b = (\tilde b_1,..., \tilde b_d)$] is $\R^k$-valued [resp. $\R^{d\times k}$-valued, resp.
$\R^{d}$-valued]. ${W}$ is a
$k$-dimensional standard Brownian motion.
They define the averaged coefficients as limits in the Ces\`aro
sense. With the additional assumption that the presumed limiting SDE has a weakly unique (in law) solution, they prove that the process
$(X^{1,x,\varepsilon}_s,\,X^{2,x,\varepsilon}_s)$ converges in
distribution towards a Markov diffusion $(X^{1,x}_s,\,X^{2,x}_s)$.
As a byproduct, they obtain an homogenization property for the
 linear PDE associated to
$(X^{1,x,\varepsilon}_s,\,X^{2,x,\varepsilon}_s)$ when the limit
Cauchy problem, associated to the limit diffusion
$(X^{1,x}_s,\,X^{2,x}_s)$, is well posed in the Sobolev space $
W_{p,\, {loc}}^{1,2}(\R_+\times\R^d)$ for each $p\geq d+2$. Here,  $
W_{p,\, {loc}}^{1,2}(\R_+\times\R^d)$ is the Sobolev space of
all functions $u(s,x)$ defined on $\R_+\times\R^d$ such that both $u$
and all the generalized derivatives $D_s u$, $D_x u$, and $D^2_{xx}
u$ belong to $L^{p}_{loc} (\R_+\times\R^d)$.

Later, the result of \cite{KK} was extended to systems of SDE-BSDE
in \cite{BEP, BDES}. Furthermore, in \cite{BEP, BDES} the
uniqueness of the averaged SDE-BSDE as well as  that of the averaged
PDE were established under appropriate conditions, building upon  the results
from \cite{K}. However, in \cite{BEP, BDES} the backward
equation does not depend on the control variable. More precisely,
the result of \cite{KK} was extended, in \cite{BEP, BDES}, to the
following SDE-BSDE.

\begin{equation}\label{SDEBSDEH1}
 \left\{
\begin{aligned}
X^{1,x,\varepsilon}_s&=x_1+\int_0^s\varphi(\frac{X^{1,x,\varepsilon}_r}
{\varepsilon},\,X^{2,x,\varepsilon}_r)dW_r,\\
X^{2,x,\varepsilon}_s&=x_2+\int_0^s \tilde b
(\frac{X^{1,x,\varepsilon}_r}{\varepsilon},\,X^{2,x,\varepsilon}_r)dr
+\int_0^s\tilde\sigma (\frac{X^{1,x,\varepsilon}_r}{\varepsilon},\,X^{2,x,\varepsilon}_r)
d {W}_r \\
Y^{t,x,\varepsilon}_s&=H(X^{x,\varepsilon}_t)+\int_s^t
f(\frac{X^{1,x,\varepsilon}_r}
{\varepsilon},\,X^{2,x,\varepsilon}_r,\,Y^{t,x,\varepsilon}_r)dr-\int_s^t
Z^{t,x,\varepsilon}_r\,dM_r^{X^{x,\varepsilon}}
\end{aligned}
\right.
\end{equation}
where $M^{X^{x,\varepsilon}}$ is the martingale part of the process
$X^{x,\varepsilon}:= (X^{1,x,\varepsilon},\,X^{2,x,\varepsilon})$.

\noindent The system of SDE-BSDE (\ref{SDEBSDEH1}) is connected to
the semilinear PDE,

\begin{equation}\label{PDEH1}
\left \{
\begin{aligned}
 \frac{\partial
v^{\varepsilon}}{\partial s}(s,\,x)&=
(\mathcal{L}^{\varepsilon} v^{\varepsilon})(t,x)
+f(\frac{x_1}{\varepsilon},\,x_2,\,v^{\varepsilon}(t,\,x)), \quad  s\geq 0
\\
v^{\varepsilon}(0,\,x)&= H(x); \quad  x = (x^1,x^2) \in
\R\times\R^{d}.
\end{aligned}
\right.
\end{equation}
where, $\mathcal{L}^{\varepsilon}$ is the infinitesimal generator
associated to the Markov process $X^{x,\varepsilon}:=
(X^{1,x,\varepsilon},\,X^{2,x,\varepsilon})$.

In the present paper we consider the situation where the coefficient
$f$ depends upon $x$, $y$ and $z$. This more general situation will force us to develop
a new methodology. That is,  the SDE-BSDE in
consideration is defined in $[0, \ t]$ by,
\begin{equation}\label{sdebsdehz0}
\left\{
\begin{aligned}
X^{1,x,\varepsilon}_s&=x_1+\int_0^s\varphi(\frac{X^{1,x,\varepsilon}_r}
{\varepsilon},\,X^{2,x,\varepsilon}_r)dW_r,\\
X^{2,x,\varepsilon}_s&=x_2+\int_0^s \tilde b
(\frac{X^{1,x,\varepsilon}_r}{\varepsilon},\,X^{2,x,\varepsilon}_r)dr
+\int_0^s\tilde\sigma (\frac{X^{1,x,\varepsilon}_r}{\varepsilon},\,X^{2,x,\varepsilon}_r)
d {W}_r \\
Y^{t,x,\varepsilon}_s&=H(X^{x,\varepsilon}_t)+\int_s^t
f(\frac{X^{1,x,\varepsilon}_r}
{\varepsilon},\,X^{2,x,\varepsilon}_r,\,Y^{t,x,\varepsilon}_r,\,Z^{t,x,\varepsilon}_r)dr-\int_s^t
Z^{t,x,\varepsilon}_r\,dM_r^{X^{x,\varepsilon}}
\end{aligned}
\right.
\end{equation}
where $M^{X^{x,\varepsilon}}$ is the martingale part of the process
$X^{x,\varepsilon}:= (X^{1,x,\varepsilon},\,X^{2,x,\varepsilon})$, i.e.
$$M^{X^{x,\varepsilon}}_s:=\int_0^s \sigma (\frac{X^{1,x,\varepsilon}_r}
{\varepsilon},\,X^{2,x,\varepsilon}_r)dW_r, \quad 0\le s\le t.$$
 If
we put for $ \ i,\,j=1,\,...,\,d$ ,\,
$$
  b:=\begin{pmatrix}0\\ \tilde b\end{pmatrix} ,\, \ \
a_{00} := \frac{1}{2}\sum_{i=1}^{k}\varphi_i^2 ,\, \ \ \tilde\sigma := (\tilde\sigma)_{ij} ,\, \ \
\sigma := \begin{pmatrix}
\varphi \\
 \tilde\sigma \\
\end{pmatrix} ,\,   \ \ \ \tilde a:=
\frac{1}{2}(\tilde\sigma\tilde\sigma^\ast) ,\, \ \ a :=
\frac{1}{2}(\sigma\sigma^\ast)
$$
(note that $a$ is a $(d+1)\times(d+1)$ matrix, whose rows and columns are indexed from $i=1$ to $i=d$, while $\tilde{a}$ is a $d\times d$ matrix),
and
$ \displaystyle
X^{x,\varepsilon}:=\begin{pmatrix}X^{1,x,\varepsilon} \\
X^{2,x,\varepsilon}\end{pmatrix},
$
 then the SDE-BSDE (\ref{sdebsdehz0}) can be rewritten in the form
\begin{equation}\label{sdebsdehz}
\left\{
\begin{aligned}
X^{x, \varepsilon}_s&= x+\int_0^s
b(\frac{X^{1,x,\varepsilon}_r}{\varepsilon},\,X^{2,x,\varepsilon}_r)dr
+\int_0^s \sigma (\frac{X^{1,x,\varepsilon}_r}{\varepsilon},\, X^{2,x,\varepsilon}_r)dW_r,\\
Y^{t,x,\varepsilon}_s&=H(X^{x,\varepsilon}_t)+\int_s^t
f(\frac{X^{1,x,\varepsilon}_r}
{\varepsilon},\,X^{2,x,\varepsilon}_r,\,Y^{t,x,\varepsilon}_r,\,Z^{t,x,\varepsilon}_r)dr
-\int_s^t Z^{t,x,\varepsilon}_r\,dM_r^{X^{x,\varepsilon}}
\end{aligned}
\right.
\end{equation}

 In this case, the nonlinear part of the PDE associated to the
SDE-BSDE (\ref{sdebsdehz}) depends on both the solution and its
gradient. More precisely, this PDE takes the form
\begin{equation}\label{pdeh}
\left \{
\begin{aligned}
 \frac{\partial
v^{\varepsilon}}{\partial s}(t,\,x)&=
(\mathcal{L}^{\varepsilon} v^{\varepsilon})(s,x)
+f(\frac{x_1}{\varepsilon},\,x_2,\,v^{\varepsilon}(s,\,x),\,
\nabla_x v^{\varepsilon}(s,\,x)), \\
v^{\varepsilon}(0,\,x)&= H(x) ,
\end{aligned}
\right.
\end{equation}
where  $\mathcal{L}^{\varepsilon}$ \ is the infinitesimal generator
associated to the Markov process  $X^{x,\varepsilon}:=
(X^{1,x,\varepsilon},\,X^{2,x,\varepsilon})$ which is more precisely
defined by


\begin{equation*}
\mathcal{L}^{\varepsilon}:=a_{00}(\frac{x_1}{\varepsilon},\,x_2)
\frac{\partial^2 }{\partial^2 x_1}+\sum_{j =1}^d
a_{0j}(\frac{x_1}{\varepsilon},\,x_2)\frac{\partial^2
}{\partial x_{1}\partial x_{2j}} +\sum_{i,\,j =1}^d
a_{ij}(\frac{x_1}{\varepsilon},\,x_2)\frac{\partial^2 }{\partial
x_{2i}\partial x_{2j}} +\sum_{i=1}^d
b_i^{(1)}(\frac{x_1}{\varepsilon},\,x_2)\frac{\partial }{\partial
x_{2i}},
\end{equation*}
$\varphi$, $\tilde\sigma $ and $\tilde b$ are the coefficients which
were defined above, $f$ and $H$ are  real valued  measurable
functions respectively defined on
$\R^{d+1}\times\R\times\R^{d+1}$ and $\R^{d+1}$.

We want to study the asymptotic behavior of the
SDE-BSDE (\ref{sdebsdehz}) when $\varepsilon \rightarrow 0$.
 Note that under suitable conditions upon the coefficients, the
function \ $\{v^{\varepsilon}(t,\,x):=Y^{\varepsilon}_0,\, t\ge0,\
x=(x_1,x_2)\in\R^{d+1}\}$ solves the PDE (\ref{pdeh}), see e. g. Remark
2.6 in \cite{EP1}. Therefore, we will also study the asymptotic
behavior of the PDE (\ref{pdeh}).

As in \cite{BEP, BDES,KK}, we consider here the averaged
coefficients as limits in the Ces\`aro sense.  Usually, the averaged coefficients are computed as means with
respected to the (unique) invariant probability measure. In our situation, due
to the fact that the fast component is null recurrent, we have no invariant probability  measure. Therefore the classical methods
do not work. Furthermore, since the variable
$Z^{\varepsilon}$ enters the generator of the backward component
and is not relatively compact in any reasonable topology, the
identification of the limit of the finite variation process of the
backward component is  rather hard to obtain. In particular the
methods used in \cite{BEP, BDES} do not work.

In order to prove that the limit problem is well posed, we
establish the existence and uniqueness for the limiting SDE-BSDE as
well as the unique solvability of the limiting PDE in the Sobolev space
$ W_{p,\, loc}^{1,2}(\R_+\times\R^d)$, \ $p\geq d+2$ .  We use
Krylov's result \cite{K} and standard arguments of BSDEs to
establish the existence and uniqueness of the limiting SDE-BSDE.  The
unique solvability of the limiting PDE is more difficult to prove. Due
to the lack of (H\"{o}lder's) regularity of the diffusion coefficient, the
pointwise estimates of the gradient can not be obtained in our
situation. To ovoid these problems, we develop a method which
consists in establishing an $L^p$-local version
 of the Calder\'{o}n-Zygmund theorem.
Our strategy is based on the $W_{p,\, loc}^{1,\,2}$--estimate for solutions of linear PDE with discontinuous coefficients proved in
\cite{DK}. We use the Gagliardo-Nirenberg interpolation inequality in order establish a $W_{p,\, loc}^{1,\,2}$-estimates  for
solution of semilinear PDEs. We then obtain a compactness
characterization of a suitable approximating sequence of PDEs from
which we derive the existence of solutions in the space
$W_{p,\, loc}^{1,\,2}$. The uniqueness is then deduced from the
uniqueness of the limiting SDE-BSDE and the Itô-Krylov formula.


 We now pass to the averaging problem.
  The lack of a reasonable compactness of $(Z^{\varepsilon})$ create some difficulties in the identification of the limits.
 Note also that, since $(Z^{\varepsilon})$  is not a semimartingale, then
 the method developed in \cite{BEP, BDES, KK} do not directly apply. To avoid these difficulties, we give an
 approach which combines PDE methods with probabilistic arguments. Indeed, building on the PDEs,
 we construct a sequence of semimartingales $(Z^{\varepsilon, n})$ that we substitute
 to $(Z^{\varepsilon})$. This allows us to use
 the method developed in  \cite{BEP, BDES, KK}. Next, we show that the problems with $(Z^{\varepsilon, n})$ and that with
 $(Z^{\varepsilon})$ average to the same limit. The limits  are obtained by combining a regularization
procedure, a stability property  and weak convergence techniques already used in \cite{BEP, BDES, D,
KK}. Let  also note that, in a periodic media, some authors
have studied the asymptotic behavior of the the PDE (\ref{pdeh}). We
refer to Gaudron and Pardoux \cite{GP} in the particular PDEs whose
nonlinearity term depends upon the gradient in a quadratic growth
manner. The case where the nonlinearity depends fully upon the
gradient have been considered by Delarue \cite{D}, who developed some of the methods
which are needed in this paper.

The paper is organized as follows: In section 2, we give the
formulation of the problem and state the main results. Sections  3 and 4 are
devoted to the proofs of the two main theorems.



\section{Formulation of the Problem and the main results}

\subsection{Notations}

For a given function $g(x)$, we define, whenever they exist, the following limits

\vskip 0.2cm
 $g^{+}(x_2):=\lim_{x_1\rightarrow
+\infty}\frac{1}{x_1}\int_0^{x_1}g(t,\,x_2)dt$,
 \ \ \ \ \
$g^{-}(x_2):=\lim_{x_1\rightarrow
-\infty}\frac{1}{x_1}\int_0^{x_1}g(t,\,x_2)dt$

 \vskip 0.15cm\noindent and  \ \
$g^{\pm}(x):=g^+(x_2)1_{\{x_1>0\}}+g^-(x_2)1_{\{x_1\leq
0\}}$.

Let \ $\displaystyle
\rho(x):=a_{00}(x)^{-1}$. The assumptions we shall make below will allow us to define the averaged
coefficients $\bar {b},\,\bar{a}$ and $\bar {f}$ by:
\begin{eqnarray}\label{defbar}
\bar {b}_i(x)&:=&\frac{(\rho
b_i)^{\pm}(x)}{\rho^{\pm}(x)},\,i=1,\,...,\,d \nonumber
\\
\bar {a}_{ij}(x)&:=&\frac{(\rho
a_{ij})^{\pm}(x)}{\rho^{\pm}(x)},
\,i,\,j=0,\,1,\,...,\,d
\\
\bar {f}(x,\,y,\,z)&:=&\frac{(\rho
f)^{\pm}(x,\,y,\,z)}{\rho^{\pm}(x)} \nonumber.
\end{eqnarray}

\vskip 0.1cm\noindent It is worth noting that $\bar{b},\,\bar{a}$ and $\bar{f}$ can be
discontinuous at $x_1=0$.

\subsection{Assumptions}
The following conditions will be used in this paper.

\vskip 0.4cm\noindent \textbf{Assumption (A) }
\begin{trivlist}
\item(A1) \ The functions $\tilde b$, $\tilde\sigma ,\,\varphi$ are uniformly Lipschitz in
$(x)$. Moreover, for each $x_1$ their derivatives in $x_2$ up to
and including second order derivatives
are bounded continuous functions of $x_2$.
\\
\item (A2) \     There exist  positive constants $\lambda$ and
$C_1$ such that for every $x$ and $\xi $, we have \
 \begin{equation*}
 \xi^*{a}\xi \geq \lambda \,\|\xi|^2 \,\
 \end{equation*}
 and
 $$ \left \{\begin{array}{l}
( i ) \  a_{00}(x)\leq C_1\\\\
  (ii
)\,\,\sum_{i=1}^{d}[\tilde a_{ii}(x)+b_i^2(x)]\leq
C_1(1+|x_2|^2)
\end{array}
\right.
$$

\vskip 0.3cm\noindent \textbf{Assumption (B) Limits in the Ces\`aro sense.}
\item (B1) \ We assume that, as $x_1$
tends to $\pm\infty$,

\vskip 0.15cm
$\displaystyle\frac{1}{x_1}\int_0^{x_1}\rho(t,\,x_2)dt $
\ (resp. $\displaystyle\frac{1}{x_1}\int_0^{x_1}D_{x_2}\rho(t,\,x_2)dt$, \ resp.
$\displaystyle\frac{1}{x_1}\int_0^{x_1}D_{x_2}^2\rho(t,\,x_2)dt$) \ tends to
\vskip 0.4cm
$\rho^{\pm}(x_2)$ \ (resp. $D_{x_2}\rho^{\pm}(x_2)$, \ resp.
$D_{x_2}^2\rho^{\pm}(x_2)$)  uniformly in $x_2$.

\vskip 0.4cm\noindent {\it We refer to
$\rho^{\pm}(x_2)$ as a limit   in the Ces\`aro sense}.

\vskip 0.3cm Here and below \
$D_{x_2}g$ and $D_{x_2}^2g$ \ respectively denote the gradient vector and
the matrix of second derivatives in $x_2$ of $g$.

\end{trivlist}

\begin{trivlist}
\item (B2)  \ For $i=0,\,...,d,\,j=1,\,...,\,d$,  the coefficients $\rho b_j$,\, $D_{x_2}(\rho b_j)$,\, $D_{x_2}^2(\rho b_j)$,\, $\rho \tilde a_{ij}$,\,
    \\ $D_{x_2}(\rho \tilde a_{ij})$,\, $D_{x_2}^2(\rho \tilde a_{ij})$
have  averages in the Ces\`aro sense.

\item (B3) \ For any function
$g\in\{\rho,\,\rho b_j,\,D_{x_2}(\rho b_j),\,D_{x_2}^2(\rho
b_j),\,\rho \tilde a_{ij},\,D_{x_2}(\rho
\tilde a_{ij},\,D_{x_2}^2(\rho \tilde a_{ij})\}$, there\\ exists a
bounded function $\alpha$ such that
\begin{eqnarray}
\left\{\begin{array}{ll}
\frac{1}{x_1}\int_0^{x_1}g(t,\,x_2)dt-g^{\pm}(x)=(1+|x_2|^2)\alpha  (x),\\\\
\label{G1} \lim_{|x_1|\longrightarrow \infty}\sup_{x_2\in
{\R}^d}|\alpha (x)|=0.
\end{array}
\right.
\end{eqnarray}
\end{trivlist}

\vfill\eject

\vskip 0.4cm\noindent \textbf{Assumption (C) }
\begin{trivlist}
\item  (C1)
There exist  $K > 0$ and $p\in \N^*$ such that
 for every $(x,y,y',z, z')\in\R^{d+1}\times\R^2\times\R^{1\times k}\times\R^{1\times k}$
\begin{eqnarray*}
\left\{\begin{array}{ll}
(i) \quad |f(x,y,z)-f(x,y',z')|\leq K(|y-y'| + |z-z'|)\\\\
 (ii)     \quad |f(x,y,z)|\leq
K(1+|x_2|^p + |y| + |z|)  \\\\
(iii) \quad   |H(x)|\leq K(1+ |x_1|^p + |x_2|^p) \
\hbox{and $H$ belongs to } \mathcal{W}_{p,\,loc}^{2}(\R^{d+1})
\end{array}
\right.
\end{eqnarray*}
\end{trivlist}

\begin{trivlist}
\item (C2)  $\rho f$ has a limit in the Ces\`aro
sense and there exists a bounded measurable  function $\beta$
such that
\begin{eqnarray}
\left\{\begin{array}{ll} \frac{1}{x_1}\int_0^{x_1}\rho
(t,\,x_2)f(t,\,x_2,\,y,\,z)dt-(\rho
f)_{\pm}(x,\,y,\,z)=(1+|x_2|^2+|y|^2+|z|^2)\beta
(x,\,y,\,z)\\\\
\label{G2} \lim_{|x_1|\rightarrow \infty}\sup_{(x_2,\,y,\,z) \in
{\R}^d\times \R\times\R^{1\times(d+1)}}|\beta
(x,\,y,\,z)|=0,
\end{array}
\right.
\end{eqnarray}
where $(\rho f)^{\pm}(x,\,y,\,z):=(\rho
f)^+(x_2,\,y,\,z)1_{\{x_1>0\}}+(\rho f)^-(x_2,\,y,\,z)1_{\{x_1\leq
0\}}$.\\

\item{(C3)}\,\, For every $x_1$, $\rho f$ has
derivatives up to  second order in $x_2,y, z$ and
these derivatives  are bounded  and  satisfy (C2).

\item{(C4)}\,\, For every $x_1$, the derivatives of $f$ in $x_2$, $y$ and $z$ up to and including
second order derivatives are bounded continuous functions.
\end{trivlist}

 Assume that \textbf{(A)}, \textbf{(B)}, \textbf{(C)} are satisfied. It is well known that:

\vskip 0.15cm For every $\varepsilon
> 0$ and every $(t,x)$, the system of SDE-BSDE (\ref{sdebsdehz}) has a unique solution
which we denote by $(X^{x, \varepsilon}_s, Y^{t, x,
\varepsilon}_s,\, Z^{t, x, \varepsilon}_s)_{0\leq s\leq t}$ such
that,

\vskip 0.15cm $\bullet$ \ $( Y^{t, x, \varepsilon},\, Z^{t, x,
\varepsilon})$ is $\mathcal{F}^{X^{x, \varepsilon}}$ adapted, where
$\mathcal{F}^{X^{x, \varepsilon}}$ denotes the filtration generated
by the process $X^{x, \varepsilon}$. More precisely, $(X^{x,
\varepsilon},\,  Y^{t, x, \varepsilon},\, Z^{t, x, \varepsilon})$ is
adapted to the filtration $\mathcal{F}^{B}$  generated by the
Brownian motion $B$.

   \vskip 0.15cm  $\bullet$ \ $\sup_{\varepsilon} \E\big(\sup_{0\leq s\leq
t}\vert Y_r^{t, x, \varepsilon}\vert^2 + \int_0^t \vert Z_r^{t, x,
 \varepsilon}\sigma(X_r)\vert^2 dr \big) < \infty $.

\vskip 0.2cm $\bullet$ \ For every $\varepsilon >0$, the  semilinear
PDE (\ref{pdeh}) has a unique  solution $v^{\varepsilon}$
in $\mathcal{C}^{1,2}$.

\vskip 0.2cm $\bullet$ \  Note that, since $a$ is uniformly
elliptic, we also have \ $\sup_{\varepsilon} \E \int_0^t \vert
Z_r^{t, x,
 \varepsilon}\vert^2 dr  < \infty $. Moreover, we have the relation \
 \begin{equation*}
v^{\varepsilon}(t,x) = Y_0^{t, x, \varepsilon}.
\end{equation*}

\vskip 0.3cm Let $\bar a$, $\bar b$ and $\bar f$ be the averaged
coefficients defined by (\ref{pdeh}). For  a fixed $(t,x)$, let $(X_s^{x},
Y_s^{t,x}, Z_s^{t,x})_{s\in[0,t]}$ denote the solution of the following system
of SDE-BSDE

\begin{eqnarray}
\label{sdebsdebar} \left\{\begin{array}{l}
X^x_s=x+\int_0^s \bar{b}(X^{x}_r)dr+\int_0^s\bar{\sigma}(X^{x}_r)dW_r,\,0\leq s\leq t.\\\\
Y^{x}_s=H(X^{x}_t)+\int_s^t
\bar{f}(X^{x}_r,\,Y^{t,x}_r,\,Z^{t,x}_r)dr-\int_s^t
Z^{t,x}_rdM^{X^{x}}_r,\,   \ \  0\leq s\leq t,
\end{array}
\right.
\end{eqnarray}
\ where $M^{X^x}$ is the martingale part of ${X^x}$.

\vskip 0.2cm The PDE associated to the averaged SDE-BSDE
(\ref{sdebsdebar}) is given by
\begin{equation}\label{pdebar}
\left\{
\begin{aligned}
\frac{\partial{ v}}{\partial s}(s,x)&=(\bar {L}v)(s,x)
+ \bar f(x,\,v(s,\,x),\,\nabla_xv(s,\,x)),\, \quad s\geq 0. \\
v(0,x)&= H(x).
\end{aligned}
\right.
\end{equation}
where  $\bar{L}$ is the infinitesimal generator associated to the
process $X^x$ and given by,
\begin{equation}\label{Lbar} \bar{L}(x):=
\sum_{i,\,j}\bar{a}_{ij}(x)\frac{\partial^2}{\partial
x_i\partial x_j}+\sum_i\bar{b}_i(x)\frac{\partial}{\partial
x_i}, \end{equation}

Our aim is show that,

\noindent 1)  equations  (\ref{sdebsdebar}) and (\ref{pdebar})  have (in some
sense) unique  solutions  $(Y^{t, x }_s,\, Z^{t,
x}_s)$ and $v$.

\noindent 2) $(X^{x, \varepsilon}_s, Y^{t, x, \varepsilon}_s,\, Z^{t, x,
\varepsilon}_s)$ converges in law to $(X^{x}_s, Y^{t, x }_s,\, Z^{t,
x}_s)$,

\noindent 3) $v^{\varepsilon}$ converges to $v$ in a topology which will be specified below.

\vskip 0.2cm According to  Khasminskii and Krylov \cite{KK} and
Krylov \cite{K}, we deduce

\begin{proposition}\label{khkr}   Assume that \textbf{(A), (B)} are satisfied. For each $x\in\R^{d+1}$, the
forward component $X^{x, \varepsilon}:=(X^{1,\, x,\,
\varepsilon},\,X^{2,\,x,\, \varepsilon})$ converges in law to the
 continuous process $X^x=(X^{1,x},\,X^{1,x})$ in
$C([0,t];\R^{d+1})$, equipped with the uniform topology. Moreover,
$X^x$ is the unique (in law) weak  solution of the forward component
of the system of equations (\ref{sdebsdebar}).
\end{proposition}

\subsection{The main results}

\begin{proposition}\label{thbsde} (Uniqueness of the averaged BSDE) \ Assume (A), (B), (C) be satisfied.
Then, for any $(t,x)\in\R_+\times\R^{d+1}$, the backward component
of the system of equations (\ref{sdebsdebar}) has a unique solution
$(Y^{t,x},\,Z^{t,x})$ such that,

(a) \ \ $(Y^{t,x},\,Z^{t,x})$ is $\mathcal{F}^X-$adapted and $
(Y_s^{t,x}, \int_s^t Z_r^{t,x}\,dM_r^{X^x})_{0\le s\le t}$ is
continuous.

(b) \ \ $\E\big( \sup_{0\le s\le t}|Y^{t,\,x}_s|^2+
\int_0^t \vert Z_r^{t,x}\sigma(X_r^x)\vert^2 dr \big)<\infty$.

(c) \ \ Moreover, \ $Y^{t,\,x}_0$ is deterministic.

  \noindent The uniqueness means that, if \ $(Y^{1}, Z^{1})$
and $(Y^{2}, Z^{2})$ are two solutions of the backward
component of \eqref{sdebsdebar} satisfying (a)--(b) then,  \
$\E\left(\sup_{0\leq s\leq t}\left|Y^{1}_s-Y^{2}_s\right|^2
+\int_0^t\left|Z_r^{1}\sigma(X_r)-Z_r^{2}\sigma(X_r)\right|^2
dr\right) = 0$
\end{proposition}

\bop  Thanks to Remark 3.5 of \cite{EP2}, it is enough to prove
existence and uniqueness of solutions for the BSDE
$$  Y^{t,\,x}_s=H(X^{x}_t)+\int_s^t
\bar f(X^{x}_r,\,Y^{t,\,x}_r,\,Z^{t,\,x}_r)dr-\int_s^t
Z^{t,\,x}_rdW_r,\, \quad 0\leq s\leq t.
$$
 Since $f$ satisfies {\bf(C)} and $\rho$ is
bounded, one can easily verify that $\bar{f}$ is uniformly Lipschitz
in $(y,\, z)$, i.e. satisfies (C1)(i). Existence and uniqueness of a solution follow
from standard results for BSDEs, see e. g. \cite{EP1}.
Finally, since $(Y_s^{t,\,x})$ is $\mathcal{F}_s^{X^x}-$adapted then
$Y^{t,\,x}_0$ is measurable with respect to a trivial
$\sigma-$algebra and hence it is deterministic. \eop

The following theorem  is closely related to the previous
proposition. It shows that the averaged PDE is uniquely solved. It will
also be used in the averaging of the SDE-BSDE as well as in the
averaging of the PDE. However, this theorem is interesting in its own
since it establishes existence, uniqueness and  $
W_{p,\text{loc}}^{1,2}([0, \ t]\times\R^d)$-regularity (for any $p\geq d+2$) of the
solution  for semilinear PDEs with discontinuous coefficients. It
extends, in some sense, the result of \cite{DK} to semilinear PDEs.

\begin{theorem}\label{thpde} Assume that (A), (B), (C) are satisfied.
Then, equation (\ref{pdebar}) has a unique  solution $v$ such
that  $v\in \ W_{p,\text{loc}}^{1,2}([0, \
t]\times\R^d)$ for any $p\geq d+2$. Moreover, this solution satisfies \ $v(t,x)
=Y_0^{t,x}$.
\end{theorem}

\vskip 0.2cm The averaging of the backward component of equation
(\ref{sdebsdehz})  is given by the following theorem.

\begin{theorem}\label{th1} [Averaging of the SDE-BSDE  (\ref{sdebsdehz})]  Assume that
(A), (B), (C) hold. Then,
 the sequence of processes $(Y_s^{t,x, \varepsilon},\,
\int_s^t Z_r^{t,x, \varepsilon}\,dM_r^{X^{\varepsilon}})_{0\le s\le
t}$ converges in law to  \ $(Y_s^{t,x}, \int_s^t
Z_r^{t,x}\,dM_r^{X^x})_{0\le s\le t}$ in $D([0,t];\R^2)$, equipped
with the ${\bf S}$--topology. Here \ $M^{X^x}$ is the martingale
part of ${X^x}$ \ and \ $(Y_s^{t,x},\, Z_s^{t,x})$ is the unique
solution of the backward component of equation (\ref{sdebsdebar}).
\end{theorem}

\begin{remark} \
In \cite{KK}, the proof is mainly based on the fact that
$X^\varepsilon$ is a semimartingale. Similarly, in \cite{BEP} the
semimartingale property which enjoy $X^\varepsilon$ and
$Y^\varepsilon$ plays an essential role, see remark 5.1 in
\cite{BEP}. If we try to follow \cite{KK} and \cite{BEP}, we need
that $Z^\varepsilon$ be a semimartingale also. Unfortunately
$Z^\varepsilon$ is not a semimartingale. Our strategy then consists
in replacing $Z^{\varepsilon}$ by an ``approximate''  semi-martingale.
 The task is to construct a continuous function $v$,
which is smooth enough such that the process
$(v(s,\,X_s),\,\nabla_xv(s,\,X_s)):=(Y_s,\,Z_s)$ is a unique
solution of the limit BSDE. To this end, by a compactness argument, we
consider the mollified coefficients $(\bar a^n,\,\bar b^n,\,\bar
f^n,\,H^n)$ and the associated solution $v^n$. Note that since our
diffusion coefficient $a$ is discontinuous, then we can not
obtain a uniform bound for $\nabla_xv^n$. We show that the sequence
$(v^n)$ can be estimated in $\mathcal W^{1,\,2}_{p,\,loc}$ uniformly
in $n$.   We then deduce a  compactness characterization of the
approximate sequence from which we derive the weak convergence
towards the function $v$. Further, we substitute $Z^{\varepsilon}$
by $\nabla_xv^n(.,\,X_.^{\varepsilon})$ in the BSDE-equation
(\ref{pdebar}).
\end{remark}



\begin{corollary}\label{th2}(Averaging of the PDE (\ref{pdeh})) Assume (A), (B), (C) hold. Then, for every $(t,x)\in\R_+\times\R^{d+1}$,
$v^\varepsilon(t,x)\to v(t,x)$, as $\varepsilon\to 0$.
\end{corollary}



\section{Proof of Theorem \ref{thpde}}
  Let
$\bar{a}_{ij}^{n},\,\bar{b}_{i}^{n},\,\bar{f}^{n},\,{H}^{n}$ denote
a regularizing sequence of $\bar{a}_{ij},\,\bar{b}_{i},\,\bar{f},\,H$ respectively. For
each $n\ge1$, $\bar{a}_{ij}^{n},\,\bar{b}_{i}^{n},\,\bar{f}^{n},\,{H}^{n}$ are
infinitely differentiable bounded functions with bounded derivatives
of every order. ${H}^{n}$ converges uniformly on compacts sets
towards $H$. Moreover   $
\bar{a}_{ij}^{n},\,\bar{b}_{i}^{n},\,\bar{f}^{n}$ converge
respectively  to $\bar{a},\,\bar{b},\,\bar{f}$ in $L^p_{loc}$ for
every $p>d+2$. We assume in addition that the assumptions (A1), (A2) and (C1)
are satisfied along the sequence, with constants which do not depend upon $n$.

Let us define
\[\displaystyle\bar{L}^{n}(x):=\sum_{i,\,j}\bar{a}_{ij}^{n}(x)\frac{\partial^2}{\partial
x_i\partial
x_j}+\sum_i\bar{b}_i^{n}(x)\frac{\partial}{\partial x_i}.\]

Consider  the sequence of PDEs on $[0, \ t]\times\R^{d+1}$,
\begin{eqnarray}
\label{PDE1} \left \{\begin{array}{l} \frac{\partial{ v^n}}{\partial
s}(s,\,x) = \bar{L}^n(x)v^n(s,x)
+\bar{f}^n(x,\,v^n(s,\,x),\,\nabla_xv^n(s,\,x))=0
\\\\
v^n(0,x)= H^n(x)
\end{array}
\right.
\end{eqnarray}
Note that, for each $n$, the PDE (\ref{PDE1}) admit a unique  solution  $v^n$ which is twice continuously differentiable in $(s, \ x)$ and three times continuously differentiable in $x$, see  e.g.
\cite{La}, Theorem 5.1, p. 320.

Using standard arguments of SDEs and
BSDEs one can show that there exists a constant $k_1$ not depending
on $n$ such that, for every $(s, \ x)$,
\begin{equation}\label{bornev^n}
 \left|v^n(s,\,x)\right|\leq k_1(1 + |x|^p).
\end{equation}
Moreover for each $n$, thanks to Theorem 7.1, chapter VII, in
Ladyzhenskaya et al. \cite{La}, or Proposition 3.3 in Ma et al.
\cite{Ma} (see also the probabilistic approach of Delarue \cite{D}
Thm. 6.1, pp. 85-89), there are constants  $k^2_n$ and $k^3_n$ such
that
\begin{equation}\label{bornegradientv^n}
 \sup_{(s,\,x)\in [0,\,t]\times
\R^{d+1}}\left|\nabla_xv^n(s,\,x)\right|\leq k_2^n\qquad \ \mbox{ and }
\qquad \sup_{(s,\,x)\in [0,\,t]\times
\R^{d+1}}\left|D^2_{xx}v^n(s,\,x)\right|\leq k_3^n
\end{equation}


\vskip 0,2cm \subsection {Compactness of the sequence $v_n$ }

\vskip 0.2cm  We now give an a priori $L^p$-bounds for the
derivatives of $v_n$.
\begin{proposition}\label{Lp-estimate}
For every $p\in[1,\,\infty[$ and $R>0$ small enough, there exists a
positive constant $C(C_1,\, K,\, p,\, R,\,t,\,k_1)$ not depending on $n$, such
that
\begin{align*}
\int_0^t \int_{
B(0,\,R/2)}\left[|\partial_sv^n|^p+|\nabla_xv^n|^p+
|D^2_{xx}v^n|^p\right]dxds\leq
C(C_1,\, K,\, p,\, R,\,t,\,k_1)
\end{align*}
\end{proposition}

    Replacing $v$ by $v-H$, the PDE (\ref{pdebar}) is reduced to a similar PDE with a null terminal datum.  Therefore, we can and do assume, throughout the proof of Proposition \ref{Lp-estimate},
 that $H = 0$.
 \par
 To establish this Proposition, we need some
preparation and lemmas.
We first recall the Gagliardo-Nirenberg interpolation inequality
which plays an important role (Theorem 3, sect. 4, Chap. 8 in Krylov
\cite{K2}, see also Theorem 7.28, Chapter VII, in Gilbarg \&
Trudinger \cite{DT}):
\begin{lemma}  \textit{(The Gagliardo-Nirenberg inequality).} \
Let $\Omega\subset\R^{d+1}$ be a bounded open set.
For any $p\ge1$, there exists a constant $ C = C(p, d, diameter(\Omega)) $  such that
for every function $\psi\in W^2_p(\Omega)$,
\begin{equation}\label{gagliardo-nirenberg}
\|\nabla_x\psi\|_{L_p(\Omega)} \leq
C\left\{\|\psi\|_{W^2_p(\Omega)}\right\}^{\frac{1}{2}}\left\{\|\psi\|_{L_p(\Omega)}\right\}^{\frac{1}{2}}.
\end{equation}
\end{lemma}

 It follows from   this inequality that,
for every $r>0$ there exists $c = c(p,\,r,\,d) > 0 $ such that for every $\varepsilon>0$,
\begin{eqnarray}
\label{TR1} \int_0^t\int_{ {B(0,\,r)}}|\nabla_x v^n(s,
x)|^pdxds &\leq &\varepsilon\,\int_0^t\int_{
{B(0,\,r)}}|D^2_{xx} v^n(s, x)|^{p}dxds\\\nonumber &+&
c(p,\,r,\,d)(1+\varepsilon^{-1})\int_0^t\int_{ {B(0,\,r)}}|
v^n(s,x)|^{p}dxds
\end{eqnarray}
Since $v^n$ is uniformly bounded on compact set, then according to the previous
inequality and the fact that $v^n$ satisfies the PDE (\ref{PDE1}), it
remains to show that for any small enough $r>0$,
 \begin{equation}\label{bornevn"}\sup_n\int_0^t\int_{
{B(0,\,r)}}|D^2_{xx} v^n(t,x)|^{p}dxdt < \infty
\end{equation}

 In order to establish the previous inequality, we  use the strategy
developed in the proof of Theorem 9.11 in Gilbarg \& Trudinger
\cite{DT}. We rewrite the PDE (\ref{PDE1}) as follows
\begin{eqnarray}\label{PDEvngn}
 \left \{\begin{array}{l} \frac{\partial{
v^{n}}}{\partial
s}(s,\,x) = \bar{a}_{ij}^{n}(x_1,\,0)\frac{\partial^2v^{n}}
{\partial x_i\partial x_j}(s,x)+g_{n}(s,\,x)=0,\,\,\,s\in(0,\,t) \\\\
v^{n}(0,x)= 0
\end{array}
\right.
\end{eqnarray}
where
\begin{eqnarray*}
g_{n}(s,\,x)&:=&\left[\bar{a}_{ij}^{n}(x)-\bar{a}_{ij}^{n}
(x_1,\,0)\right] \frac{\partial^2v^{n}}{\partial x_i\partial
x_j}(s,x)+\bar{b}_{i}^{n}(x)
\frac{\partial v^{n}}{\partial x_i}(s,x)\\
&+&\bar{f}^{n}(x,\,v^{n}(s,\,x),\,\nabla_xv^{n}(s,\,x))
\end{eqnarray*}

\noindent For $R>0$ and $s\in [0, \ t]$,   we set
\begin{itemize}

\item $Q_{s,\,t,\,R}:=[s,\,t]\times B(0,\,R)$, where $B(0,\,R)$
denotes the ball of radius $R$.

 \item $meas(Q_{s,t,R})$ denotes the Lebesgue measure of the set
$Q_{s,t,R}$.
\end{itemize}
For $\sigma\in(0,\,1)$, we put \
$\displaystyle \sigma^{\prime}:=\frac{(1+\sigma)}{2}$ \ and consider
$\eta\in\mathcal C_0^{\infty}(B_R)$ a cut--off function
$\eta:\R^{d+1}\rightarrow [0,\,1]$ satisfying the following
properties,
\begin{eqnarray*}
\left\{\begin{array}{ll}
\eta (x)&=1,\, \ \ \ \ if \ x\in B(0,\,\sigma R),\\
\eta (x)&=0,\, \ \ \ \ if \ |x|\geq \sigma' R,\\
|\nabla_x\eta (x)|&\leq {4}(1-\sigma)^{-1}R^{-1}\,\  \ \ \ \ if \ \sigma R\leq |x|\leq \sigma' R,
\\
|D^2_{xx}\eta (x)|&\leq 16(1-\sigma)^{-2}R^{-2}\,\ \ \ \ \ if \ \sigma
R\leq |x|\leq \sigma' R
\end{array}
\right.
\end{eqnarray*}
Clearly the function $u^{n}:=\eta v^{n}$ solves the PDE
\begin{eqnarray*}
 \left \{\begin{array}{l} \frac{\partial{
u^{n}}}{\partial
s}(s,\,x) = \bar{a}_{ij}^{n}(x_1,\,0)\frac{\partial^2u^{n}}{\partial x_i\partial x_j}(s,x)+G_{n}(s,\,x)=0,\,\,\,s\in(0,\,T) \\\\
u^{n}(0,x)= 0
\end{array}
\right.
\end{eqnarray*}
where, $\displaystyle
G_{n}(s,\,x):=v^{n}\bar{a}_{ij}^{n}(x_1,\,0)
\frac{\partial^2\eta}{\partial x_i\partial
x_j}+2\bar{a}_{ij}^{n}(x_1,\,0)\frac{\partial v^{n}}{\partial
x_i}\frac{\partial \eta}{\partial x_j}+\eta g_{n}(s,\,x)$

Since $\bar a^{n}$ is bounded in $x_1$ and locally Lipschitz with respect to $x_2$, uniformly w.r.t. $n$, $\bar b^{n}$
satisfies (A2) and $\bar f^{n}$ satisfies  (C1-ii), we
deduce that $G_{n}$ is bounded on $[0,t]\times\R^{d+1}$. Let $D$ be
an arbitrary bounded subset of $\R^{d+1}$. Since
$\bar{a}_{ij}^{n}(.,0)$ and $G_n$ are bounded, and $G_n$ has a
compact support, then according to Theorem 2.5  from Doyoon \& Krylov \cite{DK}, there
exists a positive constant $C= C(d,\,C_1,\, K)$ not depending on $n$ \ such that for
every $n$, we have
\begin{eqnarray}
\label{LPK} u^{n}\in {W}^{1,\,2}_{p}([0,t]\times D) \ \ \  \mbox{
and } \ \ \ \|u^{n}\|_{ W^{1,\,2}_{p}([0,t]\times D)}\leq
C\|G_{n}\|_{L^p([0,t]\times D)}.
\end{eqnarray}
From the definition of the function $\eta$, we see  that
\begin{eqnarray}
\label{TRE} \|D^2_{xx}v^{n}\|_{L_p(Q_{0,\,t,\,\sigma R})} \leq
\|D^2_{xx}u^{n}\|_{L_p(Q_{0,\,t,\,\sigma' R})}
\end{eqnarray}
According to inequalities (\ref{LPK}) and (\ref{TRE}), it remains to
estimate
 $ \displaystyle\int_0^t \int_{B(0,\,\sigma' R)}|G_{n}(s,x_1,x_2)|^pdxds $.
\par\noindent
We have
\begin{align}\label{Gn}
 \int_0^t \int_{B(0,\,\sigma' R)}|G_{n}(s,x_1,x_2)|^pdxds
\leq
 A_1+A_2+A_3
 \end{align}
 where
$$ A_1 :=
 C(p)\int_0^t\int_{B(0,\,\sigma' R)}|v^{n}|^p|\bar{a}_{ij}^{n}(x_1,\,0)|^p|
 D^2_{xx}\eta (x)|^pdxds
 $$
$$ A_2 := C(p)\int_0^t\int_{B(0,\,\sigma' R)}|\bar{a}_{ij}^{n}(x_1,\,0)|^p\,
|\nabla_x v^{n}|^p\,|\nabla_x\eta(x)|^pdxds
$$
$$
A_3 := C(p)\int_0^t\int_{B(0,\,\sigma'
R)}|g_{n}(s,\,x)|^pdxds \hskip 3.5cm
$$


   The following lemma gives estimates for $A_1, A_2$ and $A_3$.
\begin{lemma}\label{LemmaA3}
 Let \ $Q := Q_{0,t,R}$. For every $p$, there exist a positive constant $C(p)$ such that
 for every $\varepsilon>0$,
 \vskip 0.2cm\noindent
 (i) \ \ \ \ $\displaystyle
   A_1 \ \leq \
C(p)(1-\sigma)^{-2p}R^{-2p}
\int_0^t\int_{B(0,\,\sigma' R)}|v^{n}|^pdxds $ \vskip 0.2cm\noindent
(ii) $\displaystyle
 \ \ \ \   A_2 \
\leq \ C(p)(1-\sigma)^{-p}R^{-p}
\left[\varepsilon\int_0^t\int_{B(0,\,\sigma'
R)}|D^2_{xx}v^{n}|^pdxds + (1+\varepsilon^{-1})\int_0^t\int_{B(0,\,\sigma'
R)}|v^{n}|^pdxds\right] $

 \vskip 0.2cm\noindent (iii)
 $\displaystyle \ \ \ \ A_3 \
 \leq \ C(p)\left\{TR^{d+1}+R^p+(R^p+\varepsilon)
\int_0^t\int_{B(0,\,\sigma' R)}
\left|D^2_{xx}v^{n}(s,x_1,x_2)\right|^pdxds\right.$

\vskip 0.2cm\noindent $\qquad \qquad \qquad \qquad  \qquad \qquad \  \displaystyle\left.\ +\left(1+R^p\right)
(1+\varepsilon^{-1})\int_0^t\int_{B(0,\,\sigma' R)}|v^{n}|^pdxds\right\}.$

\end{lemma}
\bop $C(p)$ denotes a constant which may vary from line to line.

 \noindent Inequality $(i)$ follows from the properties of $\eta$ and the boundness of $\bar{a}_{ij}^{n}(x_1,\,0)$.


\noindent We use the properties of $\eta$, the boundedness of $\bar{a}_{ij}^{n}(x_1,\,0)$ and inequality (\ref{TR1}) to get inequality $(ii)$.
We now show inequality $(iii)$.  We have
$$\displaystyle \int_0^t\int_{B(0,\,\sigma'
R)}|g_{n}(s,\,x)|^pdxds \leq
\left(I^{n}_1+I^{n}_2+I^{n}_3\right)$$ with
\begin{eqnarray*}
I^{n}_1&:=&\int_0^t\int_{B(0,\,\sigma'
R)}\left|\bar{a}_{ij}^{n}(x)
-\bar{a}_{ij}^{n}(x_1,\,0)\right|^p
\left|\frac{\partial^2v^{n}}{\partial x_i\partial x_j}(s,x)\right|^pdxds\\\\
I^{n}_2&:=&\int_0^t\int_{B(0,\,\sigma'
R)}\left|\bar{b}_{i}^{n}(x)\right|^p
\left|\frac{\partial v^{n}}{\partial x_i}(s,x)\right|^pdxds\\\\
I^{n}_3&:=&\int_0^t\int_{B(0,\,\sigma'
R)}\left|\bar{f}^{n}(x,\,
v^{n}(s,\,x),\,\nabla_xv^{n}(s,\,x))\right|^pdxds
\end{eqnarray*}
 Since $\bar{a}_{ij}^{n}$ is uniformly Lipschitz in $x_2$, we obtain
\begin{eqnarray}
\label{LPK1} I^{n}_1&\leq
&\sup_{Q}\big(|x_2|^p\big)\int_0^t\int_{B(0,\,\sigma' R)}
\left|D^2_{xx}v^{n}(s,x_1,x_2)\right|^pdxds
\end{eqnarray}

Noticing that $\bar{b}^n$ satisfies assumption (A2-$ii$) then using inequality
(\ref{TR1}), we obtain
\begin{align}
\label{LPA2} I^{n}_2 & \leq  C_1\big(1 + \sup_{Q}|x_2|^p\big) \left[\varepsilon\int_0^t\int_{B(0,\,\sigma'
R)}|D^2_{xx}v^{n}|^pdxds\right.\\\nonumber
& \quad \ + \left.c_1\big(1+\varepsilon^{-1}\big)\int_0^t\int_{B(0,\,\sigma'
R)}|v^{n}|^pdxds\right]
\end{align}

Thanks to assumption \textbf{(C)} and inequality (\ref{TR1}) we deduce
\begin{align}
\label{LPK3} I^{n}_3&\leq
K\left(meas(Q)+\sup_{Q}(|x_2|^p)+\int_0^t\int_{B(0,\,\sigma'
R)} |v^{n}(s,x_1,x_2)|^pdxds \right.\\\nonumber
& \quad \ +\left.\varepsilon\int_0^t\int_{B(0,\,\sigma' R)}
|D^2_{xx}v^{n}(s,x_1,x_2)|^pdxds+c_1(1+\varepsilon^{-1})\int_0^t\int_{B(0,\,\sigma'
R)}|v^{n}|^pdxds\right)
\end{align}

\noindent Combining (\ref{LPK1}), (\ref{LPA2}) and (\ref{LPK3}), we
deduce the  desired result. Lemma \ref{LemmaA3} is
proved. \eop


\begin{lemma}\label{estimationv^n"}
{\bf ($L_{loc}^p$ estimate of $D^2_{xx}v^{n}$)}. \ For every
$p\in[1,\,\infty[$ and $R>0$ small enough, there exists a positive
constant $C' = C'(C_1,\, k,\,p,\,R,\,t,\,k_1)$ not depending on $n$, such that
\begin{eqnarray*}
\int_0^t\int_{B(0,\,R/2)}|D^2_{xx}v^{n}|^pdxds \leq 2R^{-2p}C'
\end{eqnarray*}
\end{lemma}
\bop
 Using inequalities (\ref{LPK}), (\ref{TRE}), (\ref{Gn}) and Lemma \ref{LemmaA3}, we show that
\begin{align*}
& (1-\sigma  )^{2p} R^{2p}  \int_0^t\int_{B(0,\,\sigma R)}|D^2_{xx}v^{n}(s,x)|^pdxds\\
&\leq  C(p)\bigg\{1 + (1-\sigma)^{p}R^{p}(1+\varepsilon^{-1}) +  (1-\sigma)^{2p}R^{2p}[1+ 2(1+\varepsilon^{-1})]
\int_0^t\int_{B(0,\,\sigma' R)}|v^{n}(s,x)|^pdxds\\
& \quad \ + (1-\sigma)^{2p}R^{2p}\bigg[\varepsilon(1-\sigma)^{-p}R^{-p} +   \sup_{Q}(|x_2|^p)(1+\varepsilon)
 + 2\varepsilon \bigg]
\int_0^t\int_{B(0,\,\sigma'
R)}\! \! \! |D^2_{xx}v^{n}(s,x)|^pdxds
\\
& \quad \ +
K(1-\sigma)^{2p}R^{2p}\big(meas(Q)+\sup_{Q}(|x_2|^p) \big)\bigg\}
\end{align*}
Using  inequality (\ref{bornev^n}) and the fact that \ $|x|\leq R$ \ in the set \ $Q:= Q_{0,t,R}$, we show that there exists a positive constant $C(C_1,\,K, R, p, k_1, \eps, meas(Q))$ such that
\begin{align*}
(1-\sigma & )^{2p}R^{2p}\int_0^t\int_{B(0,\,\sigma R)}|D^2_{xx}v^{n}|^pdxds\\
&\leq
C\left(C_1,\,K,\, R,\,p,\,k_1,\,\varepsilon,\,meas(Q)\right)\\\nonumber
& \quad \ +C(p)(1-\sigma)^{-p}R^{-p}
\left[\varepsilon(1-\sigma)^{2p}R^{2p}\int_0^t\int_{B(0,\,\sigma'
R)}|D^2_{xx}v^{n}|^pdxds\right]\\\nonumber
& \quad \ + C(p)\sup_{Q}(|x_2|^p)\left[(1-\sigma)^{2p}R^{2p}\int_0^t\int_{B(0,\,\sigma'
R)} \left|D^2_{xx}v^{n}\right|^pdxds\right]
\\\nonumber
& \quad \ + C(p)(1+\sup_{Q}| x_2|^p)
\left[\varepsilon(1-\sigma)^{2p}R^{2p}\int_0^t\int_{B(0,\sigma'
R)}|D^2_{xx}v^{n}|^pdxds\right]\\
& \quad \ + C(p)\left[\varepsilon
(1-\sigma)^{2p}R^{2p}\int_0^t\int_{B(0,\sigma'
R)}|D^2_{xx}v^{n}|^pdxds\right]
\end{align*}
Let \ $\displaystyle\bar\Lambda:= 1+ \sup_{Q}|x_2|^p$.
 We Choose $\displaystyle\varepsilon := \frac14\left\{2^{2p}\bar\Lambda C(p)\left[ (1-\sigma)^{-p}R^{-p}+2\right]\right\}^{-1}$ and $R$ be sufficiently small so that  $2^{2p}C(p)sup_{Q}(|x_2|^p)\leq \frac14$ then use the fact that $\frac{1-\sigma}{2} = 1-\sigma'$
 to obtain

\begin{align*}
(1-\sigma)^{2p}R^{2p}\int_0^t\int_{B(0,\,\sigma R)}|D^2_{xx}v^{n}|^pdxds
~ &\leq ~
\frac{1}{2}\left[(1-\sigma')^{2p}R^{2p}\int_0^t\int_{B(0,\,\sigma'
R)}|D^2_{xx}v^{n}|^pdxds\right]\\\nonumber
&  \quad \ + C(C_1,\, K,\, p,\, R,\,t,\,k_1)
\end{align*}
Passing to the \textit{sup} on $\sigma'$ and $\sigma$, we get
\begin{align*}
R^{2p}\bigg[\sup_{0<\sigma<1}(1- & \sigma)^{2p}\int_0^t\int_{B(0,\,\sigma R)}|D^2_{xx}v^{n}|^pdxds\bigg]
\\
&\leq
\frac{1}{2}R^{2p}\sup_{0<\sigma'<1}\bigg[(1-\sigma')^{2p}\int_0^t\int_{B(0,\,\sigma'
R)}|D^2_{xx}v^{n}|^pdxds\bigg]
\\
& \quad \ + C(C_1,\, K,\, p,\, R,\,t,\,k_1)
\end{align*}
It follows that
\begin{eqnarray*}
R^{2p}\left[\sup_{0<\sigma<1}(1-\sigma)^{2p}\int_0^t\int_{B(0,\,\sigma
R)}|D^2_{xx}v^{n}|^pdxds\right]
\leq 2C(C_1,\, K,\, p,\, R,\,t,\,k_1)\\
\end{eqnarray*}
The proof is finished by taking $\sigma :=1/2$. \eop

\vskip 0.15cm
\noindent\textbf{Proof of Proposition
\ref{Lp-estimate}.} Thanks to inequality (\ref{bornev^n}), inequality
(\ref{TR1}) and Lemma \ref{estimationv^n"}, we deduce that \
$sup_n\|\nabla_{x}v^n\|_{L_p([0,\,t]\times B(0,\, R/2))}$ is
 bounded. Since $v^n$ satisfies the PDE (\ref{PDE1}), we
deduce that \ $sup_n\|\partial_sv^n\|_{L_p([0,\,t]\times B(0,\,
R/2))}$ is  bounded also. Therefore, there exists a positive
constant $ C = C(C_1,\, K,\, p,\, R,\,t,\,k_1)$ such that
\begin{eqnarray}
\label{borneW2} \sup_n\int_0^t\int_{B(0,\,R/2)}\left[\vert
v^n\vert^p +
|\partial_sv^n|^p+|\nabla_xv^n|^p+|D^2_{xx}v^{n}|^p\right]dxds
\leq C
\end{eqnarray}
Proposition \ref{Lp-estimate} is proved. \eop \vskip 0,15cm

\noindent\textbf{Proof of Theorem \ref{thpde}.}  Inequalities
(\ref{borneW2}) and (\ref{bornev^n}) express that for every $R>0$ small
enough,
 $$\displaystyle \sup_{n}\|v^{n}\|_{{W}^{1,\,2}_p([0,\,t]\times B(0,\,R/2))}
 \leq C(R,\,k_1,\,t,\,p))$$
  Since, any ball $B(0,\,R^{'})$ can be covered by a finite
number of balls of radius $R/2$, and the proof of Proposition \ref{Lp-estimate}
can be easily adapted to proving the same estimate in a ball of radius $R/2$ centered around any point in
$\mathbb{R}^{d+1}$ we deduce that
\begin{equation}\label{estim1-2-p}
\sup_{n}\|v^{n}\|_{{W}^{1,\,2}_p(Q_{0,\,t,\,R'})}<\infty.
\end{equation}
Therefore $v^{n}$ converges weakly to $v$ in the space
${W}^{1,\,2}_{p}([0,\,t]\times Q)$, and $v$ solves the PDE
(\ref{pdebar}) \ $a.e$.

We now prove the uniqueness of solution in
${W}_{p,\,loc}^{1,\,2}$. Let $(X_s^{x},Y^{t,x}_s,Z^{t,x}_s)_{0\le s\le t}$ be a
solution of the FBSDE system
\begin{align}
\label{sdelim}
X^{x}_s&=x+\int_0^s \bar{b}(X^{x}_r)dr+\int_0^s\bar{\sigma}(X^{x}_r)dW_r,\, \quad 0\leq s\leq t;\\
Y^{t,x}_s&=H(X^{x}_t)+\int_s^t\bar{f}(X^{x}_r,\,Y^{t,x}_r,\,Z^{t,x}_r)dr-\int_s^t
Z^{t,x}_rdM^{X^{x}}_r,\,   \ \  0\leq s\leq t. \label{bsdelim}
\end{align}
 {\color{red} For $p\geq d+2$}, take any solution $v\in {W}^{1,\,2}_{p,\,loc}$ of the PDE  \eqref{pdebar}.   The It\^{o}-Krylov formula shows that the process
$(v(t-s,\,X_s^{x}),\,\nabla_xv(t-s,\,X_s^{x}),\,0\leq\,s\leq\,t)$ is
a solution of \eqref{bsdelim}.
Hence $v(t,x)=Y^{t,x}_0=\E(Y^{t,x}_0)$. Since \eqref{bsdelim} has a unique solution,
$v(t,x)$ is written as the expectation of a uniquely characterized functional of $(X^{x}_s)_{0\le s\le t}$. But uniqueness in law holds for \eqref{sdelim} (see Proposition \ref{khkr}), consequently
the law of $X^{x}$ is uniquely characterized,  hence the solution $v$ of \eqref{pdebar} is unique
in ${W}^{1,\,2}_{p,\,loc}$.
\eop

As consequence of Theorem \ref{thpde} and the Sobolev embedding
Theorem, we have
\begin{corollary}
\label{FU}
$v^{n}$ converges uniformly to $v$ on any compact subset of \,
$\R_+\!\times \R^{d+1}$.
\end{corollary}


\section{Proof of Theorem \ref{th1}.}

In order to simplify the notation throughout the proof of Theorem
\ref{th1}, we will suppress the superscript $x$ (resp. $(t,x)$) from the
processes $(X^{x},\, Y^{t,x},\, Z^{t,x})$ and
$(X^{x,\varepsilon},\, Y^{t,x,\varepsilon},\,
Z^{t,x,\varepsilon})$ . That is, we will respectively replace
$(X^{x},\, Y^{t,x},\, Z^{t,x})$ by $(X,\, Y,\, Z)$ and
$(X^{x,\varepsilon},\,Y^{t,x,\varepsilon},\,Z^{t,x,\varepsilon})$
by $(X^{\varepsilon},\,Y^{\varepsilon},\,Z^{\varepsilon})$.

The following lemma, can be deduced from assumption {\bf(A)}.


\begin{lemma}\label{estX}
For every   $p\ge1$ and $t>0$, there exists  constant $C(p,t)$ such
that for every $\eps>0$,
$$\E\big(\sup_{0\le s\le
t}\left[|X^{1,\eps}_s|^p+|X^{2,\eps}_s|^p + |X^{1}_s|^p+|X^{2}_s|^p\right]\big)\le C(p,t).$$
\end{lemma}



\begin{proposition}
\label{RE1} Assume that $\bf{(A),\,\,(B)}$ are satisfied. Let $\bar a$, $\bar b$, $\bar
a^n$ and $\bar b^n$ be defined as in section 3.
 Let $X=(X^1,\,X^2)$ denote the solution of the SDE
\begin{equation*}
X_s=x+\int_0^s \bar{b}(X_r)dr+\int_0^s\bar{\sigma}(X_r)dW_r,\, \ \ \
0\leq s\leq t.
\end{equation*}
Then, for every $p\geq 1$,

 (j) \ $\displaystyle\E\int_0^t|\bar{a}^n(X_r)-\bar{a}(X_r)|^pdr,\,\,\,
 \longrightarrow 0 $\,\,
 as \,\,n \,\,tends\,\, to \,\,$\infty$.

(jj) \ $\displaystyle\E\int_0^t|\bar{b}^n(X_r)-\bar{b}(X_r)|^pdr
,\,\,\, \longrightarrow 0 $\,\,as \,\,n \,\,tends\,\, to
\,\,$\infty$.
\end{proposition}
\bop
{\it Proof of $(j)$ and $(jj)$.}  \ Let $N> 0$ and put \
$D_N :=\{x\in\R^{d+1},\,|x|\leq\,N\}$. For $(g,\,g^n)\in
\{(\bar{a},\,\bar{a}^n),\,(\bar{b} ,\,\bar{b}^n)\}$, we have
\begin{eqnarray*}
\E\int_0^t|g^n(X_r)-g(X_r)|^pdr \leq  2^p \big(\E\int_0^{t}|g^n(X_r)-g(X_r)|^p \1_{\{\sup_{s\leq\, r}|X_s|\leq\, N\}}dr \\
+\E\int_0^t|g^n(X_r)-g(X_r)|^ p \1_{\{\sup_{s\leq\, r}|X_s|> N\}}dr \big)
\end{eqnarray*}
\noindent Since $\bar g$ and $g^n$ satisfy {\bf (A), (B)}, there
exists a constant $C$ which is independent of $n$ such that,
\begin{eqnarray*}
\E\int_0^t|g^n(X_r)-g(X_r)|^pdr \leq  2^p \big(\E\int_0^{t}|g^n(X_r)-g(X_r)|^p\1_{\{\sup_{s\leq\, r}|X_s|\leq\, N\}}dr \\
+\frac{C}{N^p}\E(\sup_{0\leq s\leq t}|X_s|^{2p}) \big)
\end{eqnarray*}
By Krylov's estimate, there exists a positive constant $K(t, N, d)$
which is independent of $n$ such that
\begin{eqnarray*}
\E\int_0^t|g^n(X_r)-g(X_r)|^pdr \leq K(t, N,
d+1)\|~|g^n-g|^p~\|_{L^{d+1}(D_N)}+\frac{C}{N^p}\E(\sup_{0\leq s\leq
t}|X_s|^{2p}),
\end{eqnarray*}
Passing successively to the limit in  $n$   and  $N$, we get  the desired result.  \eop





\subsubsection{Tightness of the processes $(Y^{\varepsilon},\, M^{\varepsilon}:=\int
Z^{\varepsilon}_rdM^{X^{\varepsilon}}_r)$} Recall that the process
$Y^{\varepsilon}$ is defined by
\begin{eqnarray}
\label{E9} Y^{\varepsilon}_s=H(X^{\varepsilon}_t)+\int_s^t
f(\bar{X}^{1,\,\varepsilon}_r,\,X^{2,\,\varepsilon}_r,\,Y^{\varepsilon}_r,\,Z^{\varepsilon}_r)dr-\int_s^t
Z^{\varepsilon}_r\,dM_r^{X^{\varepsilon}},
\end{eqnarray}
\noindent where
$\bar{X}^{1,\,\varepsilon}_s=\frac{X^{1,\,\varepsilon}}{\varepsilon}$.

\begin{proposition}
\label{R5} There exists a positive constant $C$ which does not
depend on $\varepsilon$ such that
\begin{eqnarray}\label{borneYepsZeps}
\sup_{\varepsilon}\left\{\E\left(\sup_{0\leq s\leq
t}\left|Y^{\varepsilon}_s\right|^2
+\int_0^t\left|Z_s^{\varepsilon}\right|^2 d\langle
M^{X^{\varepsilon}}\rangle_s\right)\right\}\leq C.
\end{eqnarray}
\end{proposition}
\bop Throughout this proof, $K$ and $C$ are positive constants which
depend only on $(s,\,t)$ and may change from line to line. According to Lemma \ref{estX} we have, for every $k\geq 1$,
\begin{eqnarray}
\label{E10} \sup_{\varepsilon}{\E}\left (\sup_{0\leq s \leq
t}\left[|X^{1,\,\varepsilon}_s|^{2k}+|X^{2,\,\varepsilon}_s|^{2k}\right]\right
) <+\infty.
\end{eqnarray}
Using It{\^o}'s formula, we get
\begin{eqnarray*}
|Y^{\varepsilon}_s|^2+\int_s^t|Z_r^{\varepsilon}|^2d{\langle\,M^{X^{\varepsilon}}\,
\rangle}_r &\leq &   |H(X^{\varepsilon}_t)|^2+K\int_s^t
|Y^{\varepsilon}_r|^2dr
+\int_s^t|f(\bar{X}^{1,\,\varepsilon}_r,\,X^{2,\,\varepsilon}_r,\,0,\,0)|^2 dr\\
&+&2C\int_s^t|Y^{\varepsilon}_r||Z^{\varepsilon}_r|dr -2\int_s^t
\langle Y^{\varepsilon}_r,\,
Z_r^{\varepsilon}dM^{X^{\varepsilon}}_s\rangle.
\end{eqnarray*}
\noindent Since
$|\sigma(\bar{X}^{1,\,\varepsilon}_r,\,X^{2,\,\varepsilon}_r)|^2
=Trace\left(\sigma\sigma^{*}
\left(\bar{X}^{1,\,\varepsilon}_r,\,X^{2,\,\varepsilon}_r\right)\right)\geq
c>0$, one has
\begin{eqnarray*}
2C|Y^{\varepsilon}_r||Z^{\varepsilon}_r| \leq
C|Y^{\varepsilon}_r|^2+\frac{1}{2}|Z^{\varepsilon}_r|^2|\sigma(\bar{X}^{1,\,\varepsilon}_r,\,X^{2,\,\varepsilon}_r)|^2.
\end{eqnarray*}
It follows that
\begin{eqnarray*}
{\E}\left (|Y^{\varepsilon}_s|^2+
\frac{1}{2}\int_s^t|Z_r^{\varepsilon}|^2d{\langle\,M^{X^{\varepsilon}}\,\rangle}_r\right)
&\leq &{\E}\left (|H(X^{\varepsilon}_t)|^2\right)+C_1{\E}\left
(\int_s^t|f(\bar{X}^{1,\,\varepsilon}_r,\,X^{2,\,\varepsilon}_r,\,0,\,0)|^2dr\right)\\
&+&K{\E} \left(\int_s^t |Y^{\varepsilon}_r|^2dr\right).
\end{eqnarray*}
According to Gronwall's Lemma, there exists a constant which does
not depend on $\varepsilon$ such that
\begin{eqnarray}
{\E}\left (|Y^{\varepsilon}_s|^2\right ) \leq C{\E}\left
(|H(X^{\varepsilon}_t)|^2
+\int_0^t|f(\bar{X}^{1,\,\varepsilon}_r,\,X^{2,\,\varepsilon}_r,\,0,\,0)|^2dr\right), \ \ \ \forall
s\in[0,\,t].\nonumber
\end{eqnarray}
We deduce that
\begin{eqnarray}
\label{E12}
 {\E}\left
(\int_s^t|Z_r^{\varepsilon}|^2d{\langle\,M^{X^{\varepsilon}}\,\rangle}_r\right
) \leq  C{\E}\left (|H(X^{\varepsilon}_t)|^2
+\int_0^t|f(\bar{X}^{1,\,\varepsilon}_r,\,X^{2,\,\varepsilon}_r,\,0,\,0)|^2dr\right)
\end{eqnarray}
Combining (\ref{E12}) and  Burkh{\"o}lder-Davis-Gundy's inequality,
we get
\begin{eqnarray*}
{\E}\left (\sup_{0\leq s \leq t}|Y^{\varepsilon}_s|^2\right ) \leq
C{\E}\left(|H(X^{\varepsilon}_t)|^2
+\int_0^t|f(\bar{X}^{1,\,\varepsilon}_r,\,X^{2,\,\varepsilon}_r,\,0,\,0)|^2dr\right).
\end{eqnarray*}
Hence,
\begin{eqnarray*}
{\E}\left (\sup_{0\leq s \leq t}|Y^{\varepsilon}_t|^2
+\frac{1}{2}\int_0^t|Z_r^{\varepsilon}|^2d{\langle\,M^{X^{\varepsilon}}\,\rangle}_r
\right) \leq  C{\E}\left
(|H(X^{\varepsilon}_t)|^2+\int_0^t|f(\bar{X}^{1,\,\varepsilon}_r,\,X^{2,\,\varepsilon}_r,\,0,\,0)|^2dr
\right)\,\,
\end{eqnarray*}
In view of  condition  (C1-$ii$ and $iii$) and inequality (\ref{E10}), the
proof is complete. \eop

\begin{proposition}
\label{R7} Let $M^{\varepsilon}_s:=\int_0^s
Z^{\varepsilon}_r\,dM_r^{X^{\varepsilon}}$. The sequence
$\left(Y^{\varepsilon},\,M^{\varepsilon}\right)_{\varepsilon>0}$ is
tight on the space
$\mathcal{D}\left([0,\,t],\,\R^L\right)\times\mathcal{D}\left([0,\,t],\,\R^{L}\right)$
 endowed with the $\bf{S}$-topology.
 \end{proposition}
\bop Since $M^{\varepsilon}$ is a martingale, then according to \cite{MZ} or \cite{J}, the Meyer-Zheng
tightness criteria is fulfilled whenever
\begin{eqnarray}
\label{ET1}
\sup_{\varepsilon}\left(CV(Y^{\varepsilon})+\E\left(\sup_{0\leq
s\leq
t}|Y^{\varepsilon}_s|+|M^{\varepsilon}_s|\right)\right)<+\infty,
\end{eqnarray}
 where $CV$ denotes the conditional variation and is defined in appendix A.

\vskip 0.1cm\noindent Clearly
\begin{eqnarray*}
CV(Y^{\varepsilon}) \leq \E\left
(\int_0^t|f(\bar{X}^{1,\,\varepsilon}_s,
\,X^{2,\,\varepsilon}_s,\,Y^{\varepsilon}_s,\,Z^{\varepsilon}_s)|ds\right).
\end{eqnarray*}
Combining condition $(C1)$ and   Proposition \ref{R5}, we
derive  (\ref{ET1}). \eop


\subsubsection{A sequence of auxiliary processes, tightness and convergence. }

\noindent For $n\in\N^*$,
 we  define a sequence of \textit{an auxiliary process}
$Z^{\varepsilon,\,n}$ by
\begin{align}\label{E7}
 Z^{\varepsilon,\,n}_s \ := \ \nabla_x v^n(t-s,\,X^{\varepsilon}_s) ,\ \ \ \ \    s\in [0, \
t]
\end{align}

\noindent We rewrite the process $Y^{\varepsilon}$ in the form,
\begin{eqnarray}
\label{E14} && Y^{\varepsilon}_s=H(X^{\varepsilon}_t) + \int_s^t
f(\bar{X}^{1,\,\varepsilon}_r,\,X^{2,\,\varepsilon}_r,\,Y^{\varepsilon}_r,\,Z^{\varepsilon,\,n}_r)dr
+   A^{\varepsilon,\,n}_t-A^{\varepsilon,\,n}_s
-(M^{\varepsilon}_t-M_s^{\varepsilon})
\end{eqnarray}
\noindent where
\begin{eqnarray}\label{Anepsilon}
\begin{array}{ll}
\displaystyle M^{\varepsilon}_s:=\int_0^s
Z^{\varepsilon}_r\,dM_r^{X^{\varepsilon}}
\\\\
\displaystyle  A^{\varepsilon,\,n}_s:=\int_0^s
\left[f(\bar{X}^{1,\,\varepsilon}_r,\,X^{2,\,\varepsilon}_r,\,Y^{\varepsilon}_r,\,Z^{\varepsilon}_r)
-f(\bar{X}^{1,\,\varepsilon}_r,\,X^{2,\,\varepsilon}_r,\,Y^{\varepsilon}_r,\,
Z^{\varepsilon,\,n}_r)\right]dr.
\end{array}
\end{eqnarray}

%
We define
\begin{eqnarray*}
\mathcal{M}^{\varepsilon,\,n}_s&:=&\int_0^s Z^{\varepsilon,\,n}_r
\sigma (\bar{X}^{1,\,\varepsilon}_r,\,X^{2,\,\varepsilon}_r)dW_r \\
&\stackrel{(\ref{E7})}=& \int_0^s \nabla_x
v^n(r,\,X^{\varepsilon}_r)
\sigma (\bar{X}^{1,\,\varepsilon}_r,\,X^{2,\,\varepsilon}_r)dW_r \\\\
\mathcal{N}^{\varepsilon,\,n}_s&:=&\int_0^s
1_{\left\{|Z^{\varepsilon}_r-Z^{\varepsilon,\,n}_r|>0\right\}}
\frac{(Z^{\varepsilon}_r-Z^{\varepsilon,\,n}_r) \sigma
(\bar{X}^{1,\,\varepsilon}_r,\,X^{2,\,\varepsilon}_r)}
{|(Z^{\varepsilon}_r-Z^{\varepsilon,\,n}_r) \sigma
(\bar{X}^{1,\,\varepsilon}_r,\,X^{2,\,\varepsilon}_r)|}dW_r \\\\
 \mathsf{L}^{\varepsilon,\,n}_s  &:=
 &  \left<\mathcal{N}^{\varepsilon,\,n},\,{M}^{\varepsilon}-
\mathcal{M}^{\varepsilon,\,n}\right>_s  \\\\
&=&\int_0^s1_{\left\{|Z^{\varepsilon}_r-Z^{\varepsilon,\,n}_r|>0\right\}}
\frac{[(Z^{\varepsilon}_r-Z^{\varepsilon,\,n}_r)
\sigma(\bar{X}^{1,\,\varepsilon}_r,\,X^{2,\,\varepsilon}_r)]
[(Z^{\varepsilon}_r-Z^{\varepsilon,\,n}_r)
\sigma(\bar{X}^{1,\,\varepsilon}_r,\,X^{2,\,\varepsilon}_r)]^{*}}
{|(Z^{\varepsilon}_r-Z^{\varepsilon,\,n}_r)\sigma
(\bar{X}^{1,\,\varepsilon}_r,\,X^{2,\,\varepsilon}_r)|}dr
\end{eqnarray*}

\begin{proposition}
\label{R6} For every $n\in\N^{*}$, the sequence
$\left(\mathcal{M}^{\varepsilon,\,n},\,
\mathcal{N}^{\varepsilon,\,n}
,\,A^{\varepsilon,\,n},\,\mathsf{L}^{\varepsilon,\,n}
\right)_{\varepsilon>0} $ is tight  on the space
$\left(\mathcal{C}\left([0,\,t],\,\R\right)\right)^4$
 endowed with the topology of uniform convergence.
\end{proposition}

 \bop  We prove the tightness of  $(\mathsf{L}^{\varepsilon,\,n})_{\varepsilon>0}$.
 Since $Z^{\varepsilon,\,n}_s \ := \ \nabla_x v^n(t-s,\,X^{\varepsilon}_s)$, then according to inequalities  (\ref{borneYepsZeps}), (\ref{bornegradientv^n}) and (\ref{E10}), we have for any $n,\, p\in\N^*$:

 \begin{equation}\label{borneXYZeps} \mbox{Max}\bigg(\sup_{\varepsilon}\E\int_0^t|Z_r^\varepsilon|^2dr  ,\, \
 \displaystyle \sup_{\varepsilon}\E\int_0^t|Z_r^{\varepsilon,\,n}|^2dr  ,\, \  \displaystyle \sup_{\varepsilon}\E\sup_{0\leq r \leq t}|X_r^{2,\varepsilon}|^pdr  \bigg)~ < ~\infty.
 \end{equation}
 We successively use  assumption \textbf{(A2)} and Schwarz's inequality to show that for any $n$
 \begin{align}\label{bornetightness}
\sup_{\varepsilon}\E\left(\sup_{|s'-s|\leq \delta}
 |\mathsf{L}_{s'}^{\varepsilon,\,n}-\mathsf{L}_s^{\varepsilon,\,n}| \right)
 & \leq
 \sup_{\varepsilon}\E\left(\sup_{|s'-s|\leq \delta}\int_s^{s'}
|(Z^{\varepsilon}_r-Z^{\varepsilon,\,n}_r)
\sigma(\bar{X}^{1,\,\varepsilon}_r,\,X^{2,\,\varepsilon}_r)|dr\right)
 \\
& \leq
 K \sup_{\varepsilon}\E\left(\sup_{r\leq t}(1+ |X_r^{2,\,\varepsilon}|)\sup_{|s'-s|\leq \delta}\int_s^{s'}
|(Z^{\varepsilon}_r-Z^{\varepsilon,\,n}_r)|dr\right) \notag
\\
& \leq
2\sqrt{\delta} K \sup_{\varepsilon}\E\bigg(\sup_{r\leq t}\big(1+ |X_r^{2,\,\varepsilon}|\big)\, \big[\int_0^t
(|Z^{\varepsilon}_r|^2+|Z^{\varepsilon,\,n}_r|^2)dr\big]^{\frac12}\bigg) \notag
\\
& \leq
 C\sqrt{\delta} . \end{align}

 Using inequality (\ref{borneXYZeps}) then letting $\delta$ tends to $0$, we deduce the tightness of  $\left(\mathsf{L}^{\varepsilon,\,n}\right)_{\varepsilon>0}$ from Theorem 7.3 in \cite{Bi}. The tightness of $(A^{\varepsilon,\,n})_{\varepsilon>0}$,  $\left(\mathcal{M}^{\varepsilon,\,n}\right)_{\varepsilon>0}$
and $\left(\mathcal{N}^{\varepsilon,\,n}\right)_{\varepsilon>0} $  can be established by similar arguments. \eop


\begin{theorem}
\label{R8}  For every $n$, there exists a continuous process
$\left(\mathcal{M}^{n},\,\mathcal{N}^{n},\,\mathsf{L}^{n},\,A^{n}\right),$
a c\`ad-l\`ag process $\left(\bar{Y},\,\bar{M} \right)$ such that
along a subsequence of $\varepsilon$, we have:
\vskip 0,10cm \noindent
$\displaystyle\left(\mathcal{M}^{\varepsilon,\,
n},\,\mathcal{N}^{\varepsilon,\,n},\, \mathsf{L}^{\varepsilon,\,n},
\,A^{\varepsilon,\,n},\,{Y}^{\varepsilon},\,{M}^{\varepsilon}
\right) \Rightarrow
\left(\mathcal{M}^{n},\,\mathcal{N}^{n},\,\mathsf{L}^{n},
\,A^{n},\,\bar{Y},\,\bar{M} \right)$ on
$\left(\mathcal{C}\left([0,\,t],\,\R\right)\right)^4\times\left(\mathcal{D}
\left([0,\,t],\,\R\right)\right)^2$
respectively endowed with the topology of the uniform convergence
and the $\bf{S}$-topology.

Moreover there exists a
countable subset $\mathsf{D}$ of $ [0,\,t]$ such that for any $k\ge1$,
$t_1,\ldots,t_k\in\mathsf{D}^c$,
\[ (Y_{t_1}^\eps,M^\eps_{t_1},\ldots,Y_{t_k}^\eps,M^\eps_{t_k})
\Rightarrow (\bar{Y}_{t_1},\bar{M}_{t_1},\ldots,\bar{Y}_{t_k},\bar{M}_{t_k}),\]
where  $\Rightarrow$ denotes the convergence in law.
\end{theorem}

\bop   \ From Propositions  \ref{R7} and \ref{R6}, the family
$\displaystyle
\left(\mathcal{M}^{\varepsilon,\,n},\,\mathcal{N}^{\varepsilon,\,n},\,\mathsf{L}^{\varepsilon,\,n},\,
A^{\varepsilon,\,n},\,{Y}^{\varepsilon},\,{M}^{\varepsilon}
\right)_{\varepsilon}$ is tight on
$\left(\mathcal{C}\left([0,\,t],\,\R\right)\right)^4\times\left(\mathcal{D}
\left([0,\,t],\,\R\right)\right)^2$, where
 the spaces are respectively endowed with the topology of the
uniform convergence and the $\bf{S}$-topology. We deduce that along
a subsequence (still denoted by $\varepsilon$), $\displaystyle
\left(\mathcal{M}^{\varepsilon,\,n},\,\mathcal{N}^{\varepsilon,\,n},\,\mathsf{L}^{\varepsilon,\,n},\,
A^{\varepsilon,\,n},\,{Y}^{\varepsilon},\,{M}^{\varepsilon}
\right)_{\varepsilon}$ converges in law on
$\left(\mathcal{C}\left([0,\,t],\,\R\right)\right)^4\times\left(\mathcal{D}
\left([0,\,t],\,\R\right)\right)^2$  to a process $\displaystyle
\left(\mathcal{M}^{n},\,\mathcal{N}^{,n},\,\mathsf{L}^{n},\,A^{n},\,\bar{Y}^n,\,\bar{M}^n\right)$.
The last statement follows from Theorem 3.1 in Jakubowski \cite{J}.
\eop

\subsubsection{The first identification of the limits in $\eps$}

 \noindent In this subsection, we will determine the equation satisfied by the limit process
 $(\bar{Y},\,\bar{M})$.
\begin{proposition}
\label{R11} Let $(\bar{Y},\,\bar{M})$, be the process defined  in
Theorem \ref{R8} as a limit (as $\varepsilon \rightarrow 0$) of
$(Y^{\varepsilon},\,M^{\varepsilon})$.  Then,\\
$(i)$ For every $ s\in[0,\,t]- \mathsf{D}$,
\begin{eqnarray}
\left\{\begin{array}{ll} \label{E15} \bar{Y}_s=H(X_t)+\int_s^t
\bar{f}(X^1_r,\,\,X^2_r,\,\bar{Y}_r,\,\nabla_x v^n(t-r,\,X_r))dr
+A^{n}_t-A^{n}_s -(\bar{M}_t -\bar{M}_s),
\\\\
 \label{E17} \E\left(\sup_{0\leq s\leq
t}|\bar{Y}_s|^2+|X^1_s|^2+|X^2_s|^2\right)\leq C.
\end{array}
\right.
\end{eqnarray}
$(ii)$ Moreover, $\bar{M}$ is $\mathcal{F}_s^n$-martingale, where
$\mathcal{F}^n_s:=\sigma{\left\{X_r,\,\bar{Y}_r,\,\bar{M}_r,\,
\mathcal{M}^n_r,\,\mathcal{N}^n_r,\,\mathsf{L}^n_r,\,A^n_r,\,0\leq
u\leq s\right\}}$ augmented with the $\P$-null sets.
\end{proposition}
To prove this proposition, we need  some lemmas. The first
one  plays a similar role to that played by the invariant measure in the periodic case. It was introduced in \cite{KK} for a forward
SDE and later adapted in \cite{BEP} to systems of SDE-BSDE in
which the generator of the backward component does not depend on the
variable $Z$. We do not provide a proof, since that of Lemma 4.7  in \cite{BEP} can be repeated word to word
(also we have a new variable).
\begin{lemma}
\label{R2} Assume {\bf (A), (B)} and (C2)-(C4). For $(x_2,\, y,\, z) \in\R^d\times\R\times\R^{d+1}$,
 let $V^{\varepsilon}(x,y,\,z)$
denote the solution of the PDE:
\begin{equation}\label{G3}
\left\{
\begin{aligned}
a_{00}(\frac{x_1}{\varepsilon},\,x_2)D^2_{x_1}u(x,y,\,z)
&=f(\frac{x_1}{\varepsilon},\,x_2,\,y,\,z)-\bar{f}(x,\,y,\,z),
\quad x_1\in\R,\\
u(0,\,x_2,y,\,z)=D_{x_1}u(0,\,x_2,y,\,z)&=0.
\end{aligned}
\right.
\end{equation}
Then, for some bounded functions $\beta_1$ and $\beta_2$ satisfying
(\ref{G2}) we have

\noindent (i)\quad $
 D_{x_1}V^{\varepsilon}(x,y,\,z)=x_1(1+|x_2|^2+|y|^2+|z|^2)\beta_1
(\frac{x_1}{\varepsilon},\,x_2,\,y,\,z)
,$\\
and the same is true with  $D_{x_1}V^{\varepsilon}$ replaced by
$D_{x_1}D_{x_2}V^{\varepsilon}$, $D_{x_1}D_yV^{\varepsilon}$ and
$D_{x_1}D_{z}V^{\varepsilon}$ \vskip 0.25cm\noindent (ii)\quad $
V^{\varepsilon}(x,y,\,z)=x_1^2(1+|x_2|^2+|y|^2+|z|^2)\beta_2
(\frac{x_1}{\varepsilon},\,x_2,\,y,\,z)
$, \\
and the same is true with $V^{\varepsilon}$ replaced by
$D_{x_2}V^{\varepsilon}$, $D_yV^{\varepsilon}$,
$D_zV^{\varepsilon}$, $D^2_{x_2}V^{\varepsilon}$,
$D^2_yV^{\varepsilon}$, $D^2_zV^{\varepsilon}$,
$D_{x_2}D_yV^{\varepsilon}$ , $D_{x_2}D_{z}V^{\varepsilon}$ and
$D_{y}D_{z}V^{\varepsilon}$.
\end{lemma}

\begin{lemma}
\label{R9} We have, for any fixed $n\geq 1$,
$$
\sup_{0\leq s\leq t}\left |\,\int_0^t\left
(f(\frac{X^{1,\,\varepsilon}_r}{\varepsilon},\,X^{2,\,\varepsilon}_r,\,
Y^{\varepsilon}_r,\nabla_x v^n(t-r,\,X^{\varepsilon}_r))
-\bar{f}(X^{1,\,\varepsilon}_r,\,X^{2,\,\varepsilon}_r,\,Y^{\varepsilon}_r,\,\nabla_x
v^n(t-r,\,X^{\varepsilon}_r))\right )dr \,\right |
$$
tends to zero  in probability  as $\varepsilon \longrightarrow 0$.
\end{lemma}
\bop We set
$$
h(\bar{X}^{1,\,\varepsilon}_s,\,X^{2,\,\varepsilon}_s,\,Y^{\varepsilon}_s,\,Z^{\varepsilon,\,n}_s)
=f(\frac{X^{1,\,\varepsilon}_s}{\varepsilon},\,X^{2,\,\varepsilon}_s,\,Y^{\varepsilon}_s,\,Z^{\varepsilon,\,n}_s)-
\bar{f}(X^{1,\,\varepsilon}_s,\,X^{2,\,\varepsilon}_s,\,Y^{\varepsilon}_s,\,Z^{\varepsilon,\,n}_s),
. $$
We shall show that  for any $0\leq s\leq t$
\begin{align*}
\lim_{\varepsilon \rightarrow 0}
\left|\int_0^sh(\bar{X}^{1,\,\varepsilon}_r,\,X^{2,\,\varepsilon}_r,\,Y^{\varepsilon}_r,\,
Z^{\varepsilon,\,n}_r)dr
\right| = 0
\end{align*}

Let $V^{\varepsilon}$ denote the solution of equation (\ref{G3}).
Note that $V^{\varepsilon}$ has  first and second derivatives in
$(x,y,z)$ which are possibly discontinuous only at $x_1=0$.
Then, as in \cite{KK},  since $\varphi^2$ is bounded away from zero,
we can use the It{\^o}-Krylov formula to get
\begin{align}\label{ETH1}
&V^{\varepsilon}(X^{1,\,\varepsilon}_s,\,X^{2,\,\varepsilon}_s,Y^\eps_s,\,Z^{\eps,\,n}_s)
=V^{\varepsilon}(x,Y^\eps_0,\,Z^{\eps,\,n}_0) \nonumber
\\
&+\int_0^s
\big[f(\frac{X^{1,\,\varepsilon}_r}{\varepsilon},\,X^{2,\,\varepsilon}_r,\,
Y^{\varepsilon}_r,\,Z^{\eps,\,n}_r)-
\bar{f}(X^{1,\,\varepsilon}_r,\,X^{2,\,
\varepsilon}_r,\,Y^{\varepsilon}_r,\,Z^{\eps,\,n}_r)\big]dr
\nonumber
\\
&+\int_0^s Trace
\big[\tilde a(\frac{X^{1,\,\varepsilon}_r}{\varepsilon},\,X^{2,\,\varepsilon}_r)
{D_{x_2}^2V^{\varepsilon}}(X^{1,\,
\varepsilon}_r,\,X^{2,\,\varepsilon}_r,Y^\eps_r,\,Z^{\eps,\,n}_r)\big]dr
\nonumber \\
&+\int_0^s [D_{x_2}V^{\varepsilon}(X^{1,\,\varepsilon}_r,X^{2,\,
\varepsilon}_r,Y^\eps_r,\,Z^{\eps,\,n}_r)\tilde b(\frac{X^{1,\,\varepsilon}_r}{\varepsilon},\,X^{2,\,
\varepsilon}_r) - D_{y}
V^\eps(X^{1,\eps}_r,X^{2,\eps}_r,Y^\eps_r)f(\frac{X^{1,\eps}_r}{\eps},
X^{2,\eps}_r,Y^\eps_r)]dr \nonumber \\
&+\int_0^s[D_{x}
V^{\varepsilon}(X^{1,\,\varepsilon}_r,\,X^{2,\,\varepsilon}_r,Y^\eps_r,\,Z^{\eps,\,n}_r)
\sigma (\frac{X^{1,\,\varepsilon}_r}{\varepsilon},\,X^{2,\,\varepsilon}_r)
 +
D_{y}
V^\eps(X^{1,\eps}_r,X^{2,\eps}_r,Y^\eps_r,\,Z^{\eps,\,n}_r)Z^\eps_r
\sigma(\frac{X^{1,\eps}_r}{\eps},X^{2,\eps}_r)]
dW_r \nonumber  \\
&+\frac{1}{2}\int_0^s D_{y}^2
V^\eps(X^{1,\eps}_r,X^{2,\eps}_r,Y^\eps_r,\,Z^{\eps,\,n}_r)Z^\eps_r
\sigma\sigma^\ast(\frac{X^{1,\eps}_r}{\eps},X^{2,\eps}_r)(Z^\eps_r)^\ast
dr \nonumber \\
&+\frac{1}{2}\int_0^s D_{x}D_{y}
V^\eps(X^{1,\eps}_r,X^{2,\eps}_r,Y^\eps_r,\,Z^{\eps,\,n}_r)
\sigma\sigma^\ast(\frac{X^{1,\eps}_r}{\eps},X^{2,\eps}_r)(Z^\eps_r)^\ast
dr\nonumber
\\
&+\frac{1}{2}\int_0^s D_{x}D_{z}
V^\eps(X^{1,\eps}_r,X^{2,\eps}_r,Y^\eps_r,\,Z^{\eps,\,n}_r)d\langle
X^{\eps},\,Z^{\eps,\,n}\rangle_r\nonumber
\\
&+\frac{1}{2}\int_0^s D_{y}D_{z}
V^\eps(X^{1,\eps}_r,X^{2,\eps}_r,Y^\eps_r,\,Z^{\eps,\,n}_r)d\langle
Y^{\eps},\,Z^{\eps,\,n}\rangle_r\nonumber
\\
&+\frac{1}{2}\int_0^s D^2_{z}
V^\eps(X^{1,\eps}_r,X^{2,\eps}_r,Y^\eps_r,\,Z^{\eps,\,n}_r)d\langle
Z^{\eps,\,n}\rangle_r\nonumber
\\
&+\int_0^s D_{z}
V^\eps(X^{1,\eps}_r,X^{2,\eps}_r,Y^\eps_r,\,Z^{\eps,\,n}_r)d\,Z^{\eps,\,n}_r
\end{align}
In view of Lemma \ref{R2} and Proposition \ref{R5},
\begin{equation*}
\lim_{\varepsilon\rightarrow 0} V^{\varepsilon}(x, Y_0^\eps,\,Z_0^{\eps,\,n})=0
\end{equation*}
Using the fact that $1 = 1_{\{|X^{1,\,\varepsilon}_s| <
\sqrt{\varepsilon}\}}+1_{\{|X^{1,\,\varepsilon}_s|\geq
\sqrt{\varepsilon}\}}$ and Lemma \ref{R2}, we obtain
\begin{align*}
\left|V^{\varepsilon}(X^{1,\,\varepsilon}_s,\,X^{2,\,\varepsilon}_s,Y_s^\eps,\,Z^{\eps,\,n}_s)\right|
&\leq
\varepsilon\left[(1+|X^{2,\,\varepsilon}_s|^2+|Y^{\varepsilon}_s|^2+|Z^{\eps,\,n}_s|^2)|
\beta_2(\frac{X^{1,\,\varepsilon}_s}
{\varepsilon},\,X^{2,\,\varepsilon}_s,\,Y^{\varepsilon}_s,\,Z^{\eps,\,n}_s)|\right]\\
&\  +1_{\{|X^{1,\,\varepsilon}_s|\geq
\sqrt{\varepsilon}\}}|X^{1,\,\varepsilon}_s|^2\left[(1+|X^{2,\,\varepsilon}_s|^2
+|Y^{\varepsilon}_s|^2+|Z^{\eps,\,n}_s|^2)|\right.\\
&\left.\beta_2(\frac{X^{1,\,\varepsilon}_s}
{\varepsilon},\,X^{2,\,\varepsilon}_s,\,Y^{\varepsilon}_s,\,Z^{\eps,\,n}_s)|\right]\\
\end{align*}
Thanks to Lemma \ref{R2} and Proposition \ref{R5}, we deduce that
$$
\E\left(\sup_{0\leq s\leq
t}|V^{\varepsilon}(X^{1,\,\varepsilon}_s,\,X^{2,\,\varepsilon}_s,
Y_s^\eps,\,Z^{\eps,\,n}_s)|\right) \leq
K\left(\varepsilon+\sup_{|x_1|\geq\sqrt{\varepsilon}}
\sup_{(x_2,\,y, z)}|\beta_2(\frac{x^1}{\varepsilon},\,x^{2},\,y,\,z)|\right)
$$
Since  $\beta_2$ satisfies (\ref{G2}), the right
hand side of the previous inequality tends to zero as
$\varepsilon\longrightarrow 0$. Similarly, one can show that each term on the lines from the third to the last one in the above identity tend to zero.
Let us detail the arguments for the term on line six, and on the term on line eight.
Let us start with the term on line 6, which is one of the most delicate ones.
\begin{align*}
&\left|\int_0^sD_{y}^2
V^\eps(X^{1,\eps}_r,X^{2,\eps}_r,Y^\eps_r,\,Z^{\eps,\,n}_r)Z^\eps_r
\sigma\sigma^\ast(\frac{X^{1,\eps}_r}{\eps},X^{2,\eps}_r)(Z^\eps_r)^\ast
dr\right|\\
&\le C\sup_{0\le r\le s}\left|D_{y}^2
V^\eps(X^{1,\eps}_r,X^{2,\eps}_r,Y^\eps_r,\,Z^{\eps,\,n}_r)\right|\text{Trace}\int_0^sZ^\eps_r
\sigma\sigma^\ast(\frac{X^{1,\eps}_r}{\eps},X^{2,\eps}_r)(Z^\eps_r)^\ast
dr
\end{align*}
Since $\{\text{Trace}\int_0^sZ^\eps_r
\sigma\sigma^\ast(\frac{X^{1,\eps}_r}{\eps},X^{2,\eps}_r)(Z^\eps_r)^\ast
dr,\ 0\le s\le t\}$ is the increasing process associated to a
martingale  which is uniformly $L^1(\P)-$integrable, its square root has a bounded expectation.
Moreover, arguing as for $V^{\varepsilon}$, one can show that
\begin{equation*}\sup_{0\le r\le s}\left|D_{y}^2
V^\eps(X^{1,\eps}_r,X^{2,\eps}_r,Y^\eps_r,\,Z^{\eps,\,n}_r)\right| \ \ \mbox{tends in probability to 0 \ as} \ \varepsilon\rightarrow 0.
\end{equation*}

We now consider the term on line 8.
Since  $\nabla_x v^n(s,\, x) \in \mathcal{C}^{1,2}$,  we use Itô's formula to get
\begin{align}\label{gradientv^n}
\nabla_x  v^n(0,\,X^{\eps}_t)=\nabla_x   v^n(t,\,X^{\eps}_0) & +
\int_0^t\Gamma(r,\,\frac{X^{1,\eps}_r}{\eps},X^{2,\eps}_r)dr \nonumber
\\
& +\int_0^tD^2_{xx}  v^n(t-r,\,X^{\eps}_r)
\sigma(\frac{X^{1,\eps}_r}{\eps},X^{2,\eps}_r)dW_r
\end{align}
where
\begin{align*} \Gamma(r,\,\frac{X^{1,\eps}_r}{\eps},X^{2,\eps}_r)&=-\partial_r
\left(\nabla_x  v^n(t-r,\,X^{\eps}_r)\right)
-D^2_{xx}  v^n(t-r,\,X^{\eps}_r)b(\frac{X^{1,\eps}_r}{\eps},X^{2,\eps}_r)\\
&+\frac{1}{2}\partial_{x,\,x,\,x}^3  v^n(t-r,\,X^{\eps}_r)\sigma\sigma^\ast
(\frac{X^{1,\eps}_r}{\eps},X^{2,\eps}_r)
\end{align*}
According to inequalities (\ref{bornev^n}) and
(\ref{bornegradientv^n}), it follows that (\ref{gradientv^n}) is
well-defined. Moreover, we have
\begin{align*}&
\frac{1}{2}\int_0^s D_{x}D_{z}
V^\eps(X^{1,\eps}_r,X^{2,\eps}_r,Y^\eps_r,\,Z^{\eps,\,n}_r)d\langle
X^{\eps},\,Z^{\eps,\,n} \rangle_r\nonumber
\\
&\leq C\sup_{0\le r\le
s}\left|D_{x}D_{z}V^\eps(X^{1,\eps}_r,X^{2,\eps}_r,Y^\eps_r,\,Z^{\eps,\,n}_r)
\right| \times \int_0^s|Trace \
\sigma\sigma^\ast(\frac{X^{1,\eps}_r}
{\eps},X^{2,\eps}_r)D^2_{xx}  v^n(t-r,\,X^{\eps}_r)|dr
\end{align*}
\noindent In view of condition {\bf(A2)}, (\ref{E10}) and the fact
that $|D^2_{xx}v^n|\leq k^n_3$, the $L^p(\P)$ norm of the
increasing process $\int_0^s|Trace \
\sigma\sigma^\ast(\frac{X^{1,\eps}_r}
{\eps},X^{2,\eps}_r)D^2_{xx}v^n(t-r,\,X^{\eps}_r)|dr$ is bounded
(by a constant not depending on $\varepsilon$), for each $p\geq 1$.
Further, the same argument as above shows that
\begin{align*}
\sup_{0\le r\le
s}\left|D_{x}D_{z}V^\eps(X^{1,\eps}_r,X^{2,\eps}_r,Y^\eps_r,\,Z^{\eps,\,n}_r)
\right|\longrightarrow \,0, \ \  \hbox{as} \ \varepsilon
\longrightarrow \,0
\end{align*}
\noindent Similarly, one can show that
\begin{align*}&
\frac{1}{2}\int_0^s D_{y}D_{z}
V^\eps(X^{1,\eps}_r,X^{2,\eps}_r,Y^\eps_r,\,Z^{\eps,\,n}_r)d\langle
Y^{\eps},\,Z^{\eps,\,n}\rangle_r\nonumber +\frac{1}{2}\int_0^s
D^2_{z}
V^\eps(X^{1,\eps}_r,X^{2,\eps}_r,Y^\eps_r,\,Z^{\eps,\,n}_r)d\langle
Z^{\eps,\,n}\rangle_r\nonumber
\\
+&\int_0^s D_{z}
V^\eps(X^{1,\eps}_r,X^{2,\eps}_r,Y^\eps_r,\,Z^{\eps,\,n}_r)d\,Z^{\eps,\,n}_r
\end{align*}
\noindent converges to zero in probability as \ $\varepsilon$ tends
to  $0$. The proof is complete. \eop

\begin{lemma}
\label{R10} For every $n\in\N^*$, the sequence of processes \
 $\displaystyle
\int_0^. \bar{f}({X^{1,\,\varepsilon}_r},\,X^{2,\,\varepsilon}_r,\,
Y^{\varepsilon}_r,\,\nabla_x  v^n(t-r,\,X^{\varepsilon}_r))dr$
  \ converges in law (as
$\varepsilon \rightarrow 0$) to the process  $\displaystyle
\int_0^.\bar{f}(X^{1}_r,\,X^{2}_r,\,\bar{Y}_r,\,\nabla_x  v^n(t-r,\,X_r))dr$
 \
on $(\mathcal{C}([0,\,t],\,\R),\,||\,||_{\infty})$.
\end{lemma}
\bop It can be performed as in \cite{BEP}-Lemma 4.9. \eop \vskip 0.2cm
\noindent {\bf Proof of Proposition \ref{R11} } Passing to the limit
in  (\ref{E14}) and using Lemma \ref{R9} and Lemma \ref{R10}, we derive assertion $(i)$. Assertion $(ii)$ can be proved by using the same argument as those of
\cite{P},  section 6. \eop

Let $\mathcal{F}^n_s := \sigma{\left\{X_r,\,\bar{Y}_r,\,\bar{M}_r,\,
\mathcal{M}^n_r,\,\mathcal{N}^n_r,\,\mathsf{L}^n_r,\,A^n_r,\quad
0\leq u\leq s\right\}}$ be the filtration generated by
\\
$(X,\,\bar{Y},\,\bar{M},\,
\mathcal{M}^n,\,\mathcal{N}^n,\,\mathsf{L}^n,\,A^n)$ and completed
by the $\P$-null sets. Combining the estimates in Proposition
\ref{R5}, inequality (\ref{E10}), Lemmas (\ref{B3}) and
(\ref{B4}) in Appendix A, we show that $\bar{M}$ is
$\mathcal{F}^n_s$-martingale.

\vskip 0,3cm  The following proposition summarizes
Proposition 6.5.2 and Corollaries 6.5.3 and 6.5.4 in
Delarue \cite{D}. We will sketch the proof for the convenience of the reader.
\begin{proposition}
\label{R12} For every $n\in\N^*$ and every  $s\in[0,\,t]$ we have
\begin{trivlist}
\item (i)$\quad [\mathcal{N}^n,\quad \bar{M}-\mathcal{M}^n]_s=\mathsf{L}^n_s$.\\
\item (ii) The process $A^n$ is of bounded variation, and,
 for every progressively measurable process \newline $\{\beta_s:\quad 0\leq s \leq t\}$ satisfying $\E\left(\int_0^t|\beta_r|^2dr\right)<+\infty$ we have for any $0\leq s \leq s' \leq t$,
\begin{align}
\label{E19} \big|\int_s^{s'}\langle \beta_r,\,dA^n_r\rangle
\big|^2
\leq  C\big(\int_s^{s'}|\beta_r|^2dr\big)\big(Trace
\big\{[\bar{M}-\int_0^.Z_r^ndM^X_r]_{s'}-[\bar{M}-
\int_0^.Z_r^ndM^X_r]_{s}\big\}\big)
\end{align}
\end{trivlist}
\end{proposition}


\bop We follow \cite{D}. Assertion $(i)$ is a consequence of Theorem \ref{R8}. We prove assertion $(ii)$. Thanks to \eqref{Anepsilon} and assumption \textbf{C}, there exists $C>0$ (which value may change from line to another) such that for every $\eps>0$, $n\in N^*$ and
$ s \leq s'\leq t$ :
\begin{equation*}
\left| A^{\varepsilon,n}_{s'} - A^{\varepsilon,\,n}_s \right| \leq C \int_s^{s'} |Z_r^{\varepsilon}-Z_r^{\varepsilon, n}|ds
\end{equation*}
Using the definitions of $M^{\eps}$,  $\mathcal{M}^{\eps, n}$, $\mathcal{N}^{\eps, n}$ and the fact that the diffusion coefficient $a$ is uniformly elliptic, we deduce that :
\begin{equation*}
\left| A^{\varepsilon,n}_{s'} - A^{\varepsilon,\,n}_s \right| \leq C \, {\rm{trace}}\big([\mathcal{N}^{\eps, n},\, M^{\eps}- \mathcal{M}^{\eps, n}]_{s'} - [\mathcal{N}^{\eps, n},\, M^{\eps}- \mathcal{M}^{\eps, n}]_{s}\big)
\end{equation*}
Using Theorem \ref{R8} and assertion $(i)$, we show that for every $n\in N^*$ and
$0\leq s \leq s'\leq t$

\begin{equation*}
\left| A^{n}_{s'} - A^{n}_s \right| \leq C \, {\rm{trace}}\big([\mathcal{N}^{ n},\, \bar{M}- \mathcal{M}^{ n}]_{s'} - [\mathcal{N}^{n},\, \bar{M} - \mathcal{M}^{ n}]_{s}\big)
\end{equation*}
 Hence, thanks to the Kunita-Watanabe inequalities, for every progressively measurable process $\beta$, satisfying \ $\E\left(\int_0^t|\beta_r|^2dr\right)<+\infty$
\begin{align*}
 \big|\int_s^{s'}\langle \beta_r,\,dA^n_r\rangle
\big|
\leq  C\big(\int_s^{s'}|\beta_r|^2 d\,{\rm{trace}}[\mathcal{N}^{n}]_r\big)^{\frac12}\big({\rm{trace}}
\big\{[\bar{M}-\int_0^.Z_r^ndM^X_r]_{s'}-[\bar{M}-
\int_0^.Z_r^ndM^X_r]_s\big\}\big)^{\frac12}
\end{align*}
Since for every $\eps>0$ and  $n\in N^*$, the process $(|\mathcal{N}^{\eps,n}|^2-s)$ is a supermartingale, then  for every  $n\in N^*$ the process $(|\mathcal{N}^{n}|^2-s)$ is also a supermartingale.
Following the proof of Theorem 4.10 of Chapter I in Kratzas \& Shreve, we deduce that \ $(|{\rm{trace}}\big([\mathcal{N}^{n}]_{s'}|- [\mathcal{N}^{n}]_s)$. This completes the proof of assertion $(ii)$. \eop

\subsubsection{Identification of the limiting BSDE in $n$ }
For $ s\in [0, \ t]$ we put
\begin{eqnarray}\label{defYnZn}
\begin{array}{ll} Y^n_s :=    v^n(t-s,\,X_s) \ \ \ \ \hbox{and}
\ \ \ \ \  Z^n_s := \nabla_x  v^n(t-s,\,X_s)
\end{array}
\end{eqnarray}
\begin{proposition}\label{1ereconv}
\label{R13}   For every  $s\in [0,\,t]- \mathsf{D}$,
\begin{equation}
 \lim_{n\rightarrow
+\infty}\left(\E\left(|Y^n_{s}-\bar{Y}_{s}|\right) +\E\left\{\left([\bar{M}-
\int_0^{.}Z^n_rdM^{X}_r]_{t}-[\bar{M}-\int_0^{.}Z^n_rdM^{X}_r]_{s}\right)\right\}\right)= 0.
\end{equation}
\end{proposition}


\vskip 0.2cm \bop
 For $R>0$, let \ $D_R:=\{x\in\R^{d+1},\,|x|\leq\,R\}$ \
and \
$\tau_R:=\inf\{r>s,\,|X_r|>R\}$, \ $inf\{\emptyset\}=\infty$.

 \vskip 0.2cm\noindent
 \textbf{Step 1:} {\it  Estimate of \
$\E\left(|Y^n_{s\wedge\tau_R}-\bar{Y}_{s\wedge\tau_R}|^2\right)$. }

 \vskip 0.2cm By Itô's formula, we have
\begin{align*}
Y^n_{s} & =  v^n(0,\,X_{t})-\int_{s}^{t}\left[\frac{\partial
 v^n}{\partial r}(t-r,\,X_r)+\bar{L}  v^n(t-r,\,X_r)\right]dr
-\int_{s}^{t}\nabla_x v^n(t-r,\,X_r)dM^{X}_r\\
& = v^n(0,\,X_{t})-\int_{s}^{t}\left[\frac{\partial  v^n}{\partial
r}(t-r,\,X_r)+\bar{L}^n  v^n(t-r,\,X_r)\right]dr\\
&\quad \ + \int_{s}^{t}\left(\bar{L}^n-\bar{L}\right)  v^n(t-r,\,X_r)dr -
\int_{s}^{t}Z_r^ndM^X_r
\end{align*}
In view of (\ref{PDE1}), (\ref{E15}) and (\ref{defYnZn}), we have
\begin{align*}
Y^n_{s}-\bar{Y}_{s}&= v^n(0,\,X_{t})-\bar {Y}_{t}
+\int_{s}^{t}\left[\,\bar{f}^n(X_r,\,Y^n_r,\,Z^n_r)-\bar{f}(X_r,\,\bar{Y}_r,\,Z^n_r)\right]dr\\
& \quad \ +\int_{s}^{t}\left(\bar{L}^n-\bar{L}\right)  v^n(t-r,\,X_r)dr
-\int_{s}^{t}dA^n_r+\int_{s}^{t}\left(d\bar{M}_r-Z^n_rdM^{X}_r\right)\\
\end{align*}
Using It\^{o}'s formula on $[s\wedge\tau_R,\,t\wedge\tau_R]$, it follows that
\begin{align}\label{Yn-Ynbar}
&\E\left(|Y^n_{s\wedge\tau_R}-\bar{Y}_{s\wedge\tau_R}|^2\right)+\E\left\{\left([\bar{M}-
\int_0^{.}Z^n_rdM^{X}_r]_{t\wedge\tau_R}-
[\bar{M}-\int_0^{.}Z^n_rdM^{X}_r]_{s\wedge\tau_R}\right)\right\}\\
&=\E\left| v^n(0,\,X_{t\wedge\tau_R})-\bar
Y_{t\wedge\tau_R}\right|^2+2\E\int_{s\wedge\tau_R}^{t\wedge\tau_R}\langle
Y^n_r-\bar{Y}_r,\quad \bar{f}^n(X_r,\,Y
^n_r,\,Z^n_r)-\bar{f}(X_r,\,\bar{Y}_r,\,Z^n_r)\rangle dr  \nonumber \\
&\quad \ +2\E\int_{s\wedge\tau_R}^{t\wedge\tau_R}\langle Y^n_r-\bar{Y}_r,\quad
\left(\bar{L}^n-\bar{L}\right)  v^n(t-r,\,X_r)\rangle
dr-2\E\int_{s\wedge\tau_R}^{t\wedge\tau_R}\langle Y^n_r-\bar{Y}_r,\quad
dA^n_r\rangle. \nonumber
\end{align}
 On one hand, since $\bar f$ is uniformly Lipschitz in the $y$-variable [thanks again to Assumption (C)-($i$)], it follows
 (where the
 constant $C$ can change from  line to line),
\begin{align}
\label {E21} 2\E\int_{s\wedge\tau_R}^{t\wedge\tau_R} & \langle
Y^n_r-\bar{Y}_r, \quad \bar{f}^n(X_r,\,Y^n_r,\,Z^n_r)
-\bar{f}(X_r,\,\bar{Y}_r,\,Z^n_r)\rangle dr
\\
&\leq ~ C \E\int_{s\wedge\tau_R}^{t\wedge\tau_R}|Y^n_r-\bar{Y}_r|^2dr
+\E\int_{s\wedge\tau_R}^{t\wedge\tau_R}|\bar{f}^n(X_r,\,{Y}^n_r,\,Z^n_r)-\bar{f}(X_r,\,{Y}^n_r,\,Z^n_r)|^2
 \nonumber
 \\
  & \leq ~ C \E\int_{s\wedge\tau_R}^{t\wedge\tau_R}|Y^n_{r\wedge\tau_R}-\bar{Y}_{r\wedge\tau_R}|^2dr
+\E\int_{0}^{t\wedge\tau_R}|\bar{f}^n(X_r,\,{Y}^n_r,\,Z^n_r)
-\bar{f}(X_r,\,{Y}^n_r,\,Z^n_r)|^2
 \nonumber
 \\
 &\leq ~ C \E\int_{s}^{t}|Y^n_{r\wedge\tau_R}-\bar{Y}_{r\wedge\tau_R}|^2dr
+\E\int_{0}^{t\wedge\tau_R}|\bar{f}^n(X_r,\,{Y}^n_r,\,Z^n_r)
-\bar{f}(X_r,\,{Y}^n_r,\,Z^n_r)|^2
\nonumber
\end{align}
 The same argument shows that
\begin{align*}
\label {E22} 2\E\int_{s\wedge\tau_R}^{t\wedge\tau_R} & \langle
Y^n_r-\bar{Y}_r,\quad \left(\bar{L}^n-\bar{L}\right)
 v^n(t-r,\,X_r)\rangle dr
 \\
&\leq 2\E\int_{s}^{t}|
Y^n_{r\wedge\tau_R}-\bar{Y}_{r\wedge\tau_R}|^2dr+\E\int_{0}^{t\wedge\tau_R}|\nabla_x v^n(t-r,\, X_r)|^2|\bar{b}^n(X_r)-\bar b(X_r)|^2dr\\
& \quad \ +\E\left(\int_{0}^{t\wedge\tau_R}|D^2_{xx} v^n(t-r,\, X_r)|^2 |\bar{a}^n(X_r)-a(X_r)|^2dr\right).
\end{align*}

For   each $n\in\N^{*}$ and $R>0$, we put
\begin{align*}
\delta_1^{n,R} &:= \E\left| v^n(t-t\wedge\tau_R,\,X_{t\wedge\tau_R})-\bar
{Y}_{t\wedge\tau_R}\right|^2+
\E\int_{s}^{t\wedge\tau_R}|\bar{f}^n(X_r,\,{Y}^n_r,\,Z^n_r)-\bar{f}(X_r,\,{Y}^n_r,\,Z^n_r)|^2 dr\\
&\quad \ +\E\int_{0}^{t\wedge\tau_R}|\nabla_x v^n(t-r,\, X_r)|^2|\bar{b}^n(X_r)-\bar b(X_r)|^2dr\\
&\quad \ +\E\left(\int_{0}^{t\wedge\tau_R}|D^2_{xx} v^n(t-r,\, X_r)|^2 |\bar{a}^n(X_r)-a(X_r)|^2dr\right).
\end{align*}

\noindent  In the other hand, we deduce from inequality (\ref{E19}),
with the choice  $\beta := Y^n-\bar Y$, that for any $\alpha>0$,
\begin{align}
2\E\left|\int_{s\wedge\tau_R}^{t\wedge\tau_R}\langle Y^n_r-\bar{Y}_r,\,dA^n_r
\rangle\right|&\leq   \frac{C}{\alpha^2}\E\left(\int_{s\wedge\tau_R}^{t\wedge\tau_R}|Y^n_r
-\bar{Y}_r|^2dr\right) \\
&\quad \ +C\alpha^2\E\left(\left\{[\bar{M}-\int_0^{.}Z^n_rdM^{X}_r]_{t\wedge\tau_R}
-[\bar{M}-\int_0^{.}
Z^n_rdM^{X}_r]_{s\wedge\tau_R}\right\}\right). \nonumber
\\
&\leq  \frac{C}{\alpha^2}\E\left(\int_{s}^{t}|Y^n_{r\wedge\tau_R}
-\bar{Y}_{r\wedge\tau_R}|^2dr\right) \nonumber
\\
&\quad \ +C\alpha^2\E\left(\left\{[\bar{M}-\int_0^{.}Z^n_rdM^{X}_r]_{t\wedge\tau_R}
-[\bar{M}-\int_0^{.}
Z^n_rdM^{X}_r]_{s\wedge\tau_R}\right\}\right). \nonumber
\end{align}
We choose $\alpha^2$ such that $C\alpha^2<\frac{1}{2}$ then we use
 identity (\ref{Yn-Ynbar}) to get
\begin{align*}
 \E\left(|Y^n_{s\wedge\tau_R}-\bar{Y}_{s\wedge\tau_R}|^2\right) &
+  \frac{1}{2}\E\left\{\left([\bar{M}-\int_0^{.}Z^n_rdM^{X}_r]_{t\wedge\tau_R}
-[\bar{M}-\int_0^{.}Z^n_rdM^{X}_r]_{s\wedge\tau_R}\right)\right\} \\
&\leq \delta_1^{n,R}+C\E\int_{s}^{t}|Y^n_{r\wedge\tau_R}-\bar{Y}_{r\wedge\tau_R}|^2dr.  \\
\end{align*}
Therefore, Gronwall's Lemma yields that
\begin{align}
\label{TEA} & \E\left(|Y^n_{s\wedge\tau_R}-\bar{Y}_{s\wedge\tau_R}|^2\right)
 +  \E\left\{\left([\bar{M}-
\int_0^{.}Z^n_rdM^{X}_r]_{t\wedge\tau_R}-
[\bar{M}-\int_0^{.}Z^n_rdM^{X}_r]_{s\wedge\tau_R}\right)\right\} \nonumber
\\
 & \hskip 1cm \leq \
K_1(C,t)\delta_1^{n,R}.
\end{align}

\vskip 0.2cm\noindent  \textbf{Step 2:} \
$\displaystyle\lim_{R\rightarrow\,+\infty}\lim_{n\rightarrow\,+\infty}\delta_1^{n,R}=0$.

\vskip 0.2cm We have \
$\delta_1^{n,R} \ = \ I_1^n+I_2^n+I^n_3$, with
\begin{align*}
I_1^n \ := \  & \E\int_{0}^{t\wedge\tau_R}|\nabla_x v^n(t-r,\, X_r)|^2|\bar{b}^n(X_r)-\bar b(X_r)|^2dr \\
 & + \E\int_{0}^{t\wedge\tau_R}|D^2_{xx} v^n(t-r,\, X_r)|^2 |\bar{a}^n(X_r)-a(X_r)|^2dr,
\end{align*}
\begin{eqnarray*}
I_2^n& \ := \ &\E\int_{0}^{t\wedge\tau_R}|\bar{f}^n(X_r,\,{Y}^n_r,\,Z^n_r)-\bar{f}(X_r,\,{Y}^n_r,\,Z^n_r)|^2 dr\\
&=&\E\int_{0}^{t\wedge\tau_R}|\bar{f}^n(X_r,\, v^n(t-r,\,X_r),\,\nabla_x v^n(t-r,\,X_r))-
\bar{f}(X_r,\, v^n(t-r,\,X_r),\,\nabla_x v^n(t-r,\,X_r))|^2 dr,\\
I^n_3& \ := \ &\E\left| v^n(t-t\wedge\tau_R,\,X_{t\wedge\tau_R})-\bar
{Y}_{t\wedge\tau_R}\right|^2.
\end{eqnarray*}
\noindent  Using H\"older's inequality, Krylov's estimate, \eqref{estim1-2-p} and  Proposition \ref{RE1}, one can show that
 $I_1^n$ tends to zero as $n$ tends to infinity.

\vskip 0.3cm\noindent  We show that  $I_2^n$ tends to 0 as $n$ tends to $\infty$. \
 Let $M>0$ and put \ $I_2^n := I_2^{n,1} +
I_2^{n,2}$, \ with
\begin{eqnarray*}
I_2^{n,1} :=
\E\int_{0}^{t\wedge\tau_R}|\bar{f}^n(X_r,\,Y^n_r,\,Z^n_r)-
\bar{f}(X_r,\,Y^n_r,\,Z^n_r)|^2 1_{\{|Y^n_r|+|Z^n_r|\leq M\}}dr
\end{eqnarray*}
and
\begin{eqnarray*}
I_2^{n,2} :=
\E\int_{0}^{t\wedge\tau_R}|\bar{f}^n(X_r,\,Y^n_r,\,Z^n_r)-
\bar{f}(X_r,\,Y^n_r,\,Z^n_r)|^2 1_{\{|Y^n_r|+|Z^n_r|> M\}}dr.
\end{eqnarray*}
 We have
\begin{eqnarray*}
I_2^{n,\,1} \leq \
\E\int_{0}^{t\wedge\tau_R}\sup_{\{|y|+|z|\leq\,M\}}|\bar{f}^n(X_r^1,\,X^2_r,\,y,\,z)
-\bar{f}(X_r^1,\,X^2_r,\,y,\,z)|^2 dr.
\end{eqnarray*}
We put  \
$h^n(x):=\sup_{\{|y|+|z|\leq\,M\}}\left|\bar{f}^n(x,\,y,\,z)-
\bar{f}(x,\,y,\,z)\right|$.

\vskip 0.15cm\noindent Thanks to
 Krylov's estimate, there exists a positive constant  $N =
 N(t,R,d)$ such that
\begin{eqnarray*}
I_2^{n,\,1}\leq\E\int_0^{t\wedge\tau_R}|h^n(X_r)|^2dr \leq
N\|h^n\|_{L^{d+2}(D_R)}^2
\end{eqnarray*}
 Since $\bar{f}^n$ and $\bar{f}$ satisfy $(C1)$,  $ Y_s^n :=  v^n(t-s,X_s)$ and $ Z_s^n := \nabla_x
v^n(t-s,X_s)$,  we get
\begin{align*}
I_2^{n,\,2} &\leq \
\E\int_{0}^{t\wedge\tau_R}(|\bar{f}^n(X_r,\,Y^n_r,\,Z^n_r)| +
|\bar{f}(X_r,\,Y^n_r,\,Z^n_r)|)^2 1_{\{|Y^n_r|+|Z^n_r|> M\}}dr
\\
&\leq \ 2K\E\int_{0}^{t\wedge\tau_R}(1 + |X_r| + |Y^n_r| +
|Z^n_r|)^2 1_{\{|Y^n_r|+|Z^n_r|> M\}}dr
\\
&\leq \ 2K\left(\E\int_{0}^{t\wedge\tau_R}(1 + |X_r| + |Y^n_r| +
|Z^n_r|)^4dr\right)^{\frac12}
\left(\E\int_{0}^{t\wedge\tau_R}1_{\{|Y^n_r|+|Z^n_r|>
M\}}dr\right)^{\frac12}
\\
&\leq \frac{2K}{M^\frac12}\left(\E\int_{0}^{t\wedge\tau_R}(1 +
|X_r|^4 + |Y^n_r|^4 + |Z^n_r|^4)dr\right)^{\frac12}
\left(\E\int_{0}^{t\wedge\tau_R}{(|Y^n_r|+|Z^n_r|)
}dr\right)^{\frac12}
\\
&\leq \frac{2K}{M^\frac12}\left(\E\int_{0}^{t\wedge\tau_R}(1 +
|X_r|^4 + |v^n(t-r, X_r)|^4 + |\nabla_x v^n(t-r,
X_r)|^4)dr\right)^{\frac12}
\\
& \hskip 1.5cm \times\left(\E\int_{0}^{t\wedge\tau_R}{(|v^n(t-r, X_r)|+|\nabla_x v^n(t-r,
X_r)|) }dr\right)^{\frac12}
\end{align*}
According to Krylov's estimate, there exists a constant $N =
N(R,t,d)$ such that
\begin{align*}
\left(\E\int_{0}^{t\wedge\tau_R}(1 + |X_r|^4 + |v^n(t-r, X_r)|^4 +
|\nabla_x v^n(t-r, X_r)|^4)dr\right)^{\frac12}
 \leq N \bigg(1  & +  R +  ||v^n||_{L^{d+2}([0, \ t]\times D_R)}^4 \\
 & + ||\nabla_x v^n||_{L^{d+2}([0, \ t]\times
 D_R)}^4 \bigg)^{\frac12}
\end{align*}
and
\begin{align*}
\left(\E\int_{0}^{t\wedge\tau_R}{(|v^n(t-r, X_r)|+|\nabla_x v^n(t-r, X_r)|)
}dr\right)^{\frac12} \leq N\bigg( & ||v^n||_{L^{d+2}([0, \ t]\times
D_R)}
\\
 & \quad + ||\nabla_x v^n||_{L^{d+2}([0, \ t]\times
 D_R)}\bigg)^{\frac12}
\end{align*}
But, thanks to \eqref{estim1-2-p}, $v^n$ and $\nabla v^n$
are  bounded  in each $L_{loc}^p([0, \ t]\times \R^{d+1})$ uniformly in $n$. We then deduce that there exists a positive constant $K_1
= K_1(t,R,d)$ such that
\begin{equation*}
\sup_n I_2^{n,\,2} \ \leq \ \frac{K_1}{M^\frac12}
\end{equation*}
Therefore,
\begin{eqnarray}\label{estimateI2n}
I_2^{n}\ \leq \ K(t,R,d)\left[\|h^n\|_{L^{d+2}(D_R)}^2 +
\frac{1}{M^\frac12}\right]
\end{eqnarray}
Passing successively to the limit in $n$ and $M$, we deduce that
$I_2^n$ tends to zero as $n$ tends to infinity.

\vskip 0.3cm  We shall show that  $I_3^n$ tends to 0 as $n$ tends to $\infty$. \ We have
\begin{eqnarray*}
I^n_3&=&\E\left| v^n(t-t\wedge\tau_R,\,X_{t\wedge\tau_R})-\bar Y_{t\wedge\tau_R}\right|^2\\
&=& \E\left| v^n(t-t\wedge\tau_R,\,X_{t\wedge\tau_R})-v(t-t\wedge\tau_R,\,X_{t\wedge\tau_R})\right|^2+
\E\left|v(t-t\wedge\tau_R,\,X_{t\wedge\tau_R})-\bar
Y_{t\wedge\tau_R}\right|^2
\end{eqnarray*}
  Since as $R$ tends to $\infty$, \ $v(t-t\wedge\tau_R,\,X_{t\wedge\tau_R})$ tends to $v(0,\,X_t) =
  H(X_t)$
   and $\bar
Y_{t\wedge\tau_R}$ tends to $\bar Y_{t} = H(X_t)$, then we pass to
the limit  first in $n$ and and next in $R$
to deduce that $I_3^n$ tends to zero as $n$ tends to infinity.
Consequently  \
$\displaystyle
\lim_{R\rightarrow\,+\infty}\lim_{n\rightarrow\,+\infty}\delta_1^{n,R}=0$.

\noindent
 Since $\tau_R$ tends increasingly to infinity as $R$ tends to infinity, then for $R$ large enough \ $t\wedge\tau_R= t$ and hence \
 $\displaystyle\lim_{n\rightarrow
+\infty}\left(\E\left(|Y^n_{s}-\bar{Y}_{s}|\right) +\E\left\{\left([\bar{M}-
\int_0^{.}Z^n_rdM^{X}_r]_{t}-[\bar{M}-\int_0^{.}Z^n_rdM^{X}_r]_{s}\right)\right\}\right)= 0$.
\eop

\vskip 0.5cm We now define
\[ Y_s:=v(t-s,X_s),\quad  Z_s:=\nabla_xv(t-s,X_s),\]
where $v$ is the solution of the PDE (\ref{pdebar}).
Note that although $\nabla_xv(\cdot,\cdot)$ is only an element of $L^p_{{loc}}([0,t]\times\R^{d+1})$ (for any $p\ge d+2$),
since $X$ is non degenerate diffusion, it follows from Krylov's estimate (see \cite{K1}) that $\nabla_xv(t-s,X_s)$
is well defined as a random element of $L^2(0,t)$.
\begin{proposition}\label{2emconv}
For every  $s\in [0,\,t]$,
 \begin{equation*}
\lim_{n\rightarrow
+\infty}\left(\E\left(|Y^n_s-Y_s|\right)+ \E\int_{s}^{t}|Z^n_s-Z_s|^2d\langle
M^X\rangle_s\right) = 0
\end{equation*}
\end{proposition}
\bop
 Since $v$ belongs to $\mathcal{W}_{p,\, loc}^{1,\, 2}$,
 then It\^o--Krylov's formula and the
 uniqueness of the backward component of equation
 (\ref{sdebsdebar}) show that for every $s\in \ [0, \ t]$,
\begin{equation}\label{Y=v}
Y_s =  v(t-s, X_s)
\end{equation}

\noindent In another hand, since
\begin{eqnarray*}
\left\{\begin{array}{l} Y_s=H(X_t)+\int_s^{t}
\bar{f}(X_r,\,Y_r,\,Z_r)dr-\int_s^t
Z_rdM^{X}_r\\\\
Y^n_{s}= v^n(0,\,X_{t})-\int_{0}^{t}\bar f^n(X_r,\, v^n(t-r,\,X_r),\,\nabla_x v^n(t-r,\,X_r))dr
+\int_{0}^{t}\left(\bar{L}^n-\bar{L}\right)  v^n(t-r,\,X_r)dr\\\\
\ \ \ \ \ \ \ \  -\int_{0}^{t}Z_r^ndM^X_r
\end{array}
\right.
\end{eqnarray*}
Using It\^{o}'s formula on $[s\wedge\tau_R,\,t\wedge\tau_R]$  then arguing as in the proof of Proposition \ref{1ereconv}, it holds that
\begin{align*}
\E|Y^n_{s\wedge\tau_R}-&  Y_{s\wedge\tau_R}|^2  +\frac{1}{2}\E\int_{s\wedge\tau_R}^{t\wedge\tau_R}|Z^n_s-Z_s|^2d\langle
M^X\rangle_s
 \\
&\leq \E\left(| v^n(t-t\wedge\tau_R,\,X_{t\wedge\tau_R})-Y_{t\wedge\tau_R}|^2\right)
\\
& \quad \ + \E\int_{s\wedge\tau_R}^{t\wedge\tau_R}\langle Y_r^n-Y_r \ , \
\bar{f}^n\big(X_r,\, v^n(t-r,\,X_r),\,v^n(t-r,\,X_r)\big)
-\bar{f}(X_r,\,Y_r,\,Z_r)\rangle ds
\\
& \quad \ + \E\int_{0}^{t\wedge\tau_R}\vert\left(\bar{L}^n-\bar{L}\right)
v^n(t-r,\,X_r)\vert^2
dr \\
& \quad \ + C\E\int_{s}^{t}|Y^n_{r\wedge\tau_R} - Y_{r\wedge\tau_R}|^2dr
\end{align*}
 Since   $(Y^n_{t},
Z^n_{t}) = (v^n(0,\,X_t),\, \nabla_x v^n(0,\,X_t))$, it follows that
\begin{align*}
\frac{}{}\E|Y^n_{s\wedge\tau_R} - & Y_{s\wedge\tau_R}   |^2+\frac{1}{2}\E\int_{s\wedge\tau_R}^{t\wedge\tau_R}|Z^n_s-Z_s|^2d\langle
M^X\rangle_s
 \\
&\leq \E\left(|
v^n(t-t\wedge\tau_R,\,X_{t\wedge\tau_R})-Y_{t\wedge\tau_R}|^2\right)
\\
&  \quad \ + \E\int_{s\wedge\tau_R}^{t\wedge\tau_R}\langle Y_r^n-Y_r \ , \ \bar{f}^n(X_r,\,
Y_r^n,\,Z_r^n) -\bar{f}(X_r,\,Y_r^n,\,Z_r^n)\rangle ds
\\
&  \quad \ + \E\int_{s\wedge\tau_R}^{t\wedge\tau_R}\langle Y_r^n-Y_r \ , \ \bar{f}(X_r,\,
Y_r^n,\,Z_r^n) -\bar{f}(X_r,\,Y_r,\,Z_r^n)\rangle ds
\\
&  \quad \ + \E\int_{s\wedge\tau_R}^{t\wedge\tau_R}\langle Y_r^n-Y_r \ , \
\bar{f}(X_r,\,Y_r,\,Z_r^n) -\bar{f}(X_r,\,Y_r,\,Z_r)\rangle ds
\\
& \quad \  + \int_{s}^{t\wedge\tau_R}\vert\left(\bar{L}^n-\bar{L}\right)
v^n(t-r,\,X_r)\vert^2
dr \\
& \quad \ + C\E\int_{s}^{t}|Y^n_{r\wedge\tau_R} - Y_{r\wedge\tau_R}|^2dr
\end{align*}
Since $\bar f$ is uniformly Lipshitz in $(y,\, z)$ with
 the same Lipshitz constants
$K$
as $f$,
then for any $\alpha
>0$ satisfying $\frac{K}{\alpha} < \frac12$, we have
\begin{align*}
\E|Y^n_{s\wedge\tau_R}- & Y_{s\wedge\tau_R}  |^2+ \left(\frac{1}{2}
-\frac{K}{\alpha}\right)\E\int_{s\wedge\tau_R}^{t\wedge\tau_R}|Z^n_r-Z_r|^2d\langle
M^X\rangle_r
\\
&\leq \E\left(|
v^n(t-t\wedge\tau_R,\,X_{t\wedge\tau_R})-Y_{t\wedge\tau_R}|^2\right)
\\
&  \quad \ + \E\int_{0}^{t\wedge\tau_R}| \bar{f}^n(X_r,\, Y_r^n,\,Z_r^n)
-\bar{f}(X_r,\,Y_r^n,\,Z_r^n)|^2 ds
\\
& \quad \ + \E\int_{0}^{t\wedge\tau_R}\vert\left(\bar{L}^n-\bar{L}\right)
v^n(t-r,\,X_r)\vert^2
dr \\
& \quad \ + ( C + K + \alpha)\E\int_{s}^{t}|Y^n_{r\wedge\tau_R} - Y_{r\wedge\tau_R}|^2dr
\end{align*}

We set
\begin{align*}
\delta_2^{n,R}  & :=\E\left(| v^n(t-t\wedge\tau_R,\,X_{t\wedge\tau_R})-Y_{t\wedge\tau_R}|^2\right)
\\
& \ \ \ \ +E\int_{0}^{t\wedge\tau_R}| \bar{f}^n(X_r,\, Y_r^n,\,Z_r^n)
-\bar{f}(X_r,\,Y_r^n,\,Z_r^n)|^2 ds \\
& \ \ \ \ +
\E\int_{0}^{t\wedge\tau_R}\vert\left(\bar{L}^n-\bar{L}\right)
v^n(t-r,\,X_r)\vert^2 dr
\end{align*}
 Arguing as for $\delta_1^{n,R}$, we show that $\displaystyle
\lim_{R\rightarrow\,+\infty}\lim_{n\rightarrow\,+\infty}\delta_2^{n,R}=0$ and the conclusion follows
as in the proof of Proposition \ref{1ereconv}. \eop

\begin{corollary}\label{cor:convY}
$\displaystyle  \P\left\{\forall s\in [0, \
t], \ \ \bar{Y}_s=v(t-s,\,X_s)\right\}=1$, which implies that  \
$(\bar{Y}_s)_{s\leq t}$ is continuous. Moreover
$\displaystyle Y^\eps\Rightarrow Y$.
\end{corollary}
\bop
Combining Propositions \ref{1ereconv} and \ref{2emconv}, we deduce that for all $s\in [0,\,t]- \mathsf{D}$,
$\bar{Y}_s=Y_s=v(s,X_s)$ a.s. Hence $\bar{Y}$ has a continuous modification, which coincides a.s. with $Y$
on $[0,t]$. But $\bar{Y}$ is c\`al\`ag, hence it is a.s. continuous and identical to $Y$.

Since $\bar{Y}$ was defined as the limit in law of an arbitrary converging subsequence of the sequence $Y^\eps$,
$\bar{Y}_s=v(s,\,X_s)$, and the law of $X$ is uniquely determined,  the law of $\{v(s,X_s),\  0\le s\le t\}$ is
uniquely determined. Consequently, the whole sequence converges : $Y^\eps\Rightarrow Y$.
\eop



\vskip 0.4cm\noindent\textbf{Proof of Corollary \ref{th2}} \ From equations (\ref{E14}) and  (\ref{E15}), we have
\begin{eqnarray*}
\left\{\begin{array}{l}
Y^{\varepsilon}_0=H(X^{\varepsilon}_t)+A^{\varepsilon\,n}_t+
\int_0^t
f(\bar{X}^{\varepsilon}_r,\,X^{2,\,\varepsilon}_r,\,Y^{\varepsilon}_r,\,Z^{\varepsilon,\,
n}_r)dr
-M^{\varepsilon}_t \\\\
\bar{Y}_0=H(X_t)+A^{n}_t+\int_0^t
\bar{f}(X_r,\,\bar{Y}_r,\,Z^{n}_r)dr-\bar{M}_t
\end{array}
\right.
\end{eqnarray*}
By Corollary \ref{cor:convY} and the continuity of the projection at the final time $t\not\in D$ : $y\mapsto y_t$, we deduce from the above two identitites that
$Y^{\varepsilon}_0$ converges towards $\bar{Y}_0$ in distribution.
Moreover, since $Y^{\varepsilon}_0,\, \bar{Y}_0$
 are deterministic, we deduce that \ $\lim_{\varepsilon \rightarrow
0}Y_0^{\varepsilon}=\bar{Y}_0=Y_0$. That is, by using the non simplified
notation,
\[ Y^{t,x,\eps}_0\to Y^{t,x}_0. \]
In other words, as $\eps\to0$, \[v^{\varepsilon}(t,\,x)\to v(t,\,x).\]
  \eop



\appendix \section{Appendix: S-topology}
The $\bf{S}$-topology has been introduced by Jakubowski \cite{J}
 as a topology defined on the Skorohod space of c\`adl\`ag
functions: $\mathcal{D}([0,\,T];\,\R)$. This topology is weaker than
the Skorohod topology but tightness criteria are easier to
establish. These criteria are the same as the one used in
Meyer-Zheng   \cite{MZ}. \\     Let $N^{a,\,b}(z)$
denotes the number of up-crossing of the function
$z\in\mathcal{D}([0,\,T];\,\R)$ in a given level $a<b$. We recall
some facts about the $\bf{S}$-topology.
\begin{proposition}(A criteria for S-tight). A sequence
$(Y^{\varepsilon})_{\varepsilon>0}$ is S-tight if and only
if it is relatively compact on the S-topology.\\
Let $(Y^{\varepsilon})_{\varepsilon>0}$ be a family of stochastic
processes in $\mathcal{D}([0,\,T];\,\R)$. Then this family is tight
for the S-topology if and only if
$(\|Y^{\varepsilon}\|_{\infty})_{\varepsilon>0}$ and
$(N^{a,\,b}(Y^{\varepsilon}))_{\varepsilon>0}$ are tight for each
$a<b$.
\end{proposition}
\noindent Let
$\left(\Omega,\,\mathcal{F},\,\P,\,(\mathcal{F}_t)_{t\geq 0}\right)$
be a stochastic basis. If $(Y)_{0\leq t\leq T}$ is a process in
$\mathcal{D}([0,\,T];\,\R)$ such that $Y_t$ is integrable for any
$t$, {\it the  conditional variation of $Y$} is defined by
$$
CV(Y)=\sup_{0\leq t_1<...<t_n=T,\,partition\, of\, [0,\,T]}
\sum_{i=1}^{n-1}\E[|\E[Y_{t_{i+1}}-Y_{t_i}\left|\right.\mathcal{F}_{t_i}]|].
$$
The process is call $quasimartingale$ if $CV(Y)<+\infty$. When $Y$
is a $\mathcal{F}_t$-martingale, $ CV(Y)=0$. A variation of Doob
inequality (cf. lemma 3, p.359 in Meyer and Zheng \cite{MZ}, where it is
assumed that $Y_T=0$) implies that
$$
\P\left[\sup_{t\in[0,\,T]}|Y_t|\geq k\right]\leq\frac{2}{k}\left
(CV(Y)+\E\left[\sup_{t\in[0,\,T]}|Y_t|\right] \right),
$$
$$
\E\left[N^{a,\,b}(Y)\right]\leq\frac{1}{b-a}\left
(|a|+CV(Y)+\E\left[\sup_{t\in[0,\,T]}|Y_t|\right]\right ).
$$
It follows that a sequence $(Y^{\varepsilon})_{\varepsilon>0}$ is
S-tight if
$$
\sup_{\varepsilon>0}\left(CV(Y^{\varepsilon})+\E\left[\sup_{t\in[0,\,T]}|Y_t^{\varepsilon}|\right]\right)<+\infty.
$$
\begin{theorem}
\label{B1} Let $(Y^{\varepsilon})_{\varepsilon>0}$ be a S-tight
family of stochastic process in $\mathcal{D}([0,\,T];\,\R)$. Then
there exists a sequence $(\varepsilon_k)_{k\in\N}$ decreasing to
zero, some process $Y\in \mathcal{D}([0,\,T];\,\R)$ and a countable
subset $D\in [0,\,T]$ such that for any $n$ and any
$(t_1,\,...,\,t_n)\in[0,\,T]\backslash D$,
$$
(Y^{\varepsilon_k}_{t_1},\,...,\,Y^{\varepsilon_k}_{t_n})\stackrel{\mathcal{D}ist}{\longrightarrow}
(Y_{t_1},\,...,\,Y_{t_n})
$$
\end{theorem}
\begin{remark}
\label{B2} The projection
:$\pi_T\,\,y\in(\mathcal{D}([0,\,T];\,\R),\,S)\mapsto\,y(T)$is
 continuous (see Remark 2.4, p.8 in Jakubowski,1997), but $y\mapsto\,y(t)$ is not continuous for each $0\leq t\leq T$.
\end{remark}
\begin{lemma}
\label{B3} Let $(Y^{\varepsilon},\,M^{\varepsilon})$ be a
multidimensional process in
$\mathcal{D}([0,\,T];\,\R^{p})\,(p\in\N^{*})$ converging to
$(Y,\,M)$ in the S-topology. Let
$(\mathcal{F}_t^{X^{\varepsilon}})_{t\geq 0}$ (resp.
$(\mathcal{F}_t^{X})_{t\geq 0}$) be the minimal complete admissible
filtration for $X^{\varepsilon}$ (resp.$X$). We assume that
$\sup_{\varepsilon>0}\E\left[\sup_{0\leq t\leq
T}|M^{\varepsilon}_t|^2\right]<C_T\,\,\forall T>0,\,M^{\varepsilon}$
is a $\mathcal{F}^{X^{\varepsilon}}$-martingale and $M$ is a
$\mathcal{F}^{X}$-adapted. Then $M$ is a
$\mathcal{F}^{X}$-martingale.
\end{lemma}
\begin{lemma}
\label{B4} Let $(Y^{\varepsilon})_{\varepsilon>0}$ be a sequence of
process converging weakly in $\mathcal{D}([0,\,T];\,\R^p)$  to $Y$.
We assume that $\sup_{\varepsilon>0}\E\left[\sup_{0\leq t\leq
T}|Y^{\varepsilon}_t|^2\right]<+\infty$. Hence, for any $t\geq
0,\,E\left[\sup_{0\leq t\leq T}|Y_t|^2\right]<+\infty$.
\end{lemma}


\begin{thebibliography}{99}


\bibitem{BEP}Bahlali, K., Elouaflin, A.,  Pardoux, E.  Homogenization of semilinear PDEs
with discontinuous averaged coefficients. \textit{EJP, Vol 14
(2009), paper no. 18, pages 477-499}.

\bibitem{BDES}Bahlali, K., Elouaflin, A., M. A. Diop, A. Said. A singular perturbation for non-divergence form  semilinear PDEs with
discontinuous effective coefficients. \textit{Preprint} (2012).

\bibitem{B} Bahlali, K. Existence and uniqueness of solutions for BSDEs
with locally Lipschitz coefficient.  Electron. Comm. Probab.  7
(2002), 169--179

 \bibitem{BP} Benchérif-Madani, A.; Pardoux, É. Homogenization of a
semilinear parabolic PDE with locally periodic coefficients: a
probabilistic approach.  ESAIM Probab. Stat.  11, 385--411
(electronic),  2007.

\bibitem{BLP} Bensoussan, A.; Lions, J.-L.; Papanicolaou, G. Asymptotic
analysis for periodic structures. Studies in Mathematics and Its
Applications, 5. North-Holland, Amsterdam-New York, 1978.

\bibitem{Bi} Billingsley, P. Convergence of probability measures, 2nd ed.,
Wiley, 1999.

\bibitem{BI} Buckdahn, R.; Ichihara, N. Limit theorem for controlled
backward SDEs and homogenization of Hamilton-Jacobi-Bellman
equations.  Appl. Math. Optim.  51,  no. 1, 1--33, 2005.

\bibitem{BHP} Buckdahn, R.; Hu, Y.; Peng, S. Probabilistic approach to
homogenization of viscosity solutions of parabolic PDEs.  NoDEA
Nonlinear Differential Equations Appl.  6,  no. 4, 395--411, 1999.

\bibitem{BH} Buckdahn, R.; Hu, Y. Probabilistic approach to homogenizations
of systems of quasilinear parabolic PDEs with periodic structures.
Nonlinear Anal.  32,  no. 5, 609--619, 1998.


\bibitem{CCKS} Caffarelli, L., Crandall, M.G., Kocan, M.,
\'{S}wiech, A. On viscosity solutions of fully nonlinear equations
with measurable ingredients. \textit{Comm. Pure Appl. Math.} 49,
365-397, 1996.

\bibitem{CKLS}  Crandall, M.G., Kocan, M., Lions, P. L.,
\'{S}wiech, A. Existence results for boundary problems for uniformly
elliptic and parabolic fully nonlinear equations. \textit{Electronic
Journal of Differential equations.}, No. 1-20, 1999.

\bibitem{D} Delarue, F.  Equations diff\'erentielles stochastiques progressives-r\'etrogrades.
 Application \`a l'homogen\'eisation des EDP quasi-lin\'eaires.
 \textit{Th\`ese de Doctorat}, (2002), Aix Marseille Université. Formerly Université de Provence, Aix-Marseille I.

\bibitem{Darticle} Delarue, F. Auxiliary SDEs for homogenization of quasilinear
PDEs with periodic coefficients.  Ann. Probab.  32  (2004),  no. 3B,
2305--2361.
\bibitem{DK} Doyoom, K., Krylov, N. Parabolic equation with measurable coefficients. \textit{Potential analysis}, 26, (2006), 345-361.
\bibitem{GP} Gaudron, G., Pardoux, E.   EDSR, convergence en loi et homogen\'eisation d'EDP paraboliques s\'emi-lin\'eaires. \textit{To appear in Anna. Inst. H. Poincar\'e}, (2001)
\bibitem{DT} Gilbarg, D., Trudinger, N.S. (1983): Elliptic partial differential equations of second order, {\it Second edition. Grundlehren der Mathematischen Wissenschaften}, 224, Springer-Verlag, Berlin.
\

\bibitem{EK} El Karoui, N. Backward stochastic differential equations a
general introduction, in Backward stochastic differential equations,
N. El Karoui and L. Mazliak Edts, {\it Pitman Research Notes in
Mathematics Series 364}, 7-27, 1997.

\bibitem{EO} Essaky, E. H.; Ouknine, Y. Averaging of backward stochastic
differential equations and homogenization of partial differential
equations with periodic coefficients.  Stoch. Anal. Appl.  24,  no.
2, 277--301, 2006.

\bibitem{FR}  Freidlin M. Functional integration
and partial differential equations. {\sl Annals of Mathematics
Studies}, 109, Princeton University Press, Princeton, 1985.

\bibitem{I} Ichihara, N. A stochastic representation for fully nonlinear
PDEs and its application to homogenization.  J. Math. Sci. Univ.
Tokyo  12,  no. 3, 467--492, 2005.

\bibitem{J} Jakubowski, A. A non-Skorohod topology on the Skorohod space.
\textit{Electron. J. Probab.} 2 , paper no. 4, pp.1-21, 1997.

\bibitem{JKO} Jikov, V. V.; Kozlov, S. M.; Ole\u\i nik, O. A.
Homogenization of differential operators and integral functionals.
Translated from the Russian by G. A. Yosifian. Springer, Berlin,
1994.

\bibitem{KK} Khasminskii, R; Krylov, N. V. On averaging principle for
diffusion processes with null-recurrent fast component.
\textit{Stochastic Processes and their applications}, 93, 229-240,
2001.

\bibitem{K1} Krylov, N. V.  Controlled Diffusion Processes, (A. B. Aries,
translator), Applications of Mathematics, Vol. 14, Springer-Verlag,
New York Berlin, 1980.

\bibitem{K} Krylov, N. V. On weak uniqueness for some diffusions with
discontinuous coefficients. \textit{Stochastic Processes and their
applications}, 113, 37-64, 2004.

\bibitem{K2} Krylov, N. V. Lectures on Elliptic and Parabolic PDEs in Sobolev Space.
Graduate Studies in Mathematics, 96. American Mathematical Society,
Providence, RI, 2008.

\bibitem{La} Ladyzhenskaya O.A., Solonnikov V.A., Ural'tseva N.N. Linera and quasi-linear
Equations of Parabolic type. \textit{ American Mathematical Society,
Providence, RI,} 1968.

\bibitem{L} Lejay, A. A probabilistic approach to the homogenization of
divergence-form operators in periodic media.  \textit{Asymptot.
Anal.}  28
 no. 2, 151--162, 2001.
\bibitem{Ma} Ma, J., Protter, P., Yong, J. Solving forward backward stochastic differential equations explicitly: a four step scheme. \textit{Probab. Theory Related Fields} 98,(1994), 339-359.

\bibitem{MZ} Meyer, P. A., Zheng, W. A. Tightness criteria for laws of
semimartingales. \textit{Ann. Inst. H. Poincar\'e Probab. Statist}.
20, (4), 217-248, 1984.

\bibitem{PAN} Pankov, A. G--convergence and homogenization of nonlinear
partial differential operators. Mathematics and Its Applications,
422. Kluwer, Dordrecht, 1997.

\bibitem{EP1} Pardoux, E. Backward stochastic differential equations and
viscosity solutions of systems of semilinear parabolic and elliptic
PDEs of second order, in {\it Stochastic analysis and related topics
VI, Geilo 1996}, Progr. Probab., vol. {\bf42}, Birkh\"auser, Boston,
MA, pp. 79--127, 1998

\bibitem{EP2}  Pardoux, E.  BSDEs, weak convergence and homogenization
of semilinear PDEs in {\it F. H Clarke and R. J. Stern (eds.),
Nonlinear Analysis, Differential Equations and Control, }503-549{\it
. Kluwer Academic Publishers.}, 1999.

\bibitem{P} Pardoux, E. Homogenization of linear and semilinear second
order parabolic PDEs with periodic coefficients: A probabilistic
approach. \textit{Journal of Functional Analysis} 167,  498-520,
1999.

\bibitem{PV} Pardoux, E., Veretennikov, A.Y, Averaging of backward SDEs
with application to semi-linear PDEs. \textit{Stochastic and
Stochastic Rep.}. 60,  255-270, 1999.


\end{thebibliography}
\end{document}